\documentclass[a4paper,12pt]{article}
\usepackage{amssymb}
\usepackage{amsthm}
\usepackage{graphics}
\usepackage[centertags]{amsmath}
\usepackage{amsfonts}
\usepackage{amssymb}
\usepackage{amsthm}
\usepackage{graphics}
\usepackage{newlfont}

\usepackage{amsmath}
\usepackage{color}

\newcommand{\eps}{\varepsilon}

\newcommand {\de}{\Delta}

\newcommand {\ga}{\gamma}

\def\<{\langle}
\def\>{\rangle}

\def\om{\omega}
\def\p{\phi}

\newcommand{\rt}{\widetilde{R}}

\newcommand{\oHs}{{\kern.2ex \overline{\kern-0.2ex H^*}}}

\newcommand{\mha}{multiplier Hopf algbera }
\newcommand{\deh}{\hat{\Delta }}
\newcommand{\Ah}{\hat{A}}
\newcommand{\q}{q}

\newcommand{\eg}{\varepsilon_G}
\newcommand{\et}{\tilde{\varepsilon}}

\newcommand{\mt}{\tilde{\mu}}
\newcommand{\st}{\tilde{S}}
\newcommand{\ttt}{\tilde{T}}
\newcommand{\at}{\tilde{A}}

\newcommand{\dt}{\tilde{\Delta}}

\newcommand{\del}{\delta }

\newcommand{\ag}{A_{G} }
\newcommand{\rg}{R_{G} }

\newcommand{\mhgc}{multiplier Hopf $G-$coalgebra }
\newcommand{\complex}{\mathbb{C}}
\newcommand{\dg}{\Delta_{G} }
\newcommand{\sg}{\sigma_G}

\newcommand{\me}{\mathcal{D} }
\newcommand{\pg}{p\in G }
\newcommand{\fa}{\varphi }
\newcommand{\fag}{\varphi_G }
\newcommand{\fah}{$\varphi=\{\varphi_{p}\}_{\pg}$ }
\newcommand{\pah}{$\psi=\{\psi_{p}\}_{\pg}$ }
\newcommand{\fpq}{\varphi_{pq}}
\newcommand{\pa}{\psi}
\newcommand{\ppq}{\psi_{pq}}
\newcommand{\Apq}{A_{pq}}
\newcommand{\dpq}{\Delta_{p,q}}

\def\de{\Delta}

\theoremstyle{plain}
\newtheorem{thm}{Theorem}
\newtheorem{cor}{Corollary}[section]
\newtheorem{lem}[cor]{Lemma}
\newtheorem{prop}[cor]{Proposition}
\theoremstyle{definition}
\newtheorem{defn}[cor]{Definition}

\theoremstyle{remark}
\newtheorem{rem}{Remark}

\begin{document}
\title{\textbf{Multiplier Hopf group coalgebras from algebraic and analytical point of views}}
\author{A.Hegazi\\
Mathematics department, Faculty of Science, Mansoura University,
Egypt\\
hegazi@mans.edu.eg\\
A.Elhafz\\
Mathematics department, Faculty of Education, Suez Canal University,
Egypt\\
a.elhafz@hotmail.com}
\date{}
\maketitle
\begin{abstract}
The Multiplier Hopf Group Coalgebra was introduced by Hegazi in 2002 \cite{10} as a generalization of Hope group caolgebra, introduced by Turaev in 2000 \cite{7}, in the non-unital case. We prove that the concepts introduced by A.Van Daele  in constructing multiplier Hopf algebra \cite{4} can be adapted to serve again in our construction.
A multiplier Hopf group coalgebra is a family of algebras $A=\{A_{\alpha}\}_{\alpha \in \pi}$, ($\pi$ is a discrete group) equipped with a family of homomorphisms $\Delta=\{\Delta_{\alpha,\beta}:A_{\alpha\beta}\longrightarrow M(A_{\alpha}\otimes A_{\beta})\}_{\alpha,\beta \in \pi}$  which is called a comultiplication under some conditions, where $M(A_{\alpha}\otimes A_{\beta})$ is the multiplier algebra of $A_{\alpha}\otimes A_{\beta}$.\\
 In 2003 A. Van Daele suggest a new approach to study the same structure by consider the direct sum of the algebras $A_p$'s which will be a multiplier Hopf algebra called later group cograded multiplier Hope algebra \cite{11}. And hence there exist a one to one correspondence between multiplier Hopf Group Coalgebra and group cograded multiplier Hopf algebra. By using this one-one correspondence  we  studied multiplier Hopf Group Coalgebra    \\

\end{abstract}
\section{Introduction}

In 1992, Alfons Van Daele \cite{4,5} introduced the notion of multiplier Hopf algebra which is considered a generalization of a usual Hopf algebra. The theory of multiplier Hopf algebra provided us with a mathematical tool for studying noncompact quantum groups with Haar measure.\\
 Recently, Quasitriangular Hopf $\pi-$coalgebras were introduced by Turaev \cite{7}. He  showed that they give rise to crossed $\pi-$categories.Viralizier \cite{8} studied the algebraic properties of the Hopf $\pi-$coalgebras. He showed the existence of integrals and trace for such coalgebras and  generalized the main properties of the quasitriangular Hopf algebras to the setting of Hopf $\pi-$coalgebras.\\
Now let us give some  basic definitions which we will need in our
work:
Assume that A is an associative algebra over $\complex$. By a left multiplier of A we mean a linear map $ L:A \longrightarrow A$ such that $L(ab)=L(a) b$ for all $ a,b \in A$. Similarly, by a right multiplier we mean a linear map $R:A\longrightarrow A$ such that $R(ab)=aR(b)$ for all $a,b\in A$. A multiplier (double multiplier)  on $A$ is a pair $(L,R) $ such that $L$ is a left multiplier, $R$ is a right multiplier and $R(a)b=aL(b)$ for all $a,b\in A$ .\\
\indent we denoted by $L(A),R(A)$, and $M(A)$ that the set of all left, right, and multipliers of $A$. It is clear that the  composition of maps makes these vector spaces into algebras. \\
\indent If the product of $A$ is nondegenerate (i.e. if $ab=0$ for all $a \in A$ then $b=0$ and if  $ab=0$ for all $b \in A$ then $a=0$) then the maps $a\longrightarrow \cdot a$, $a\longrightarrow a\cdot$ and $a\longrightarrow (\cdot a,a\cdot)$ give a natural imbedding of $A$ into $L(A),R(A)$ and $M(A)$ respectively.\\
 \indent Since the tensor product of two nondegenerate algebras is again nondegenerate then $A\otimes A \subseteq M(A)\otimes M(A) \subseteq M(A\otimes A)$.\\
\indent If a homomorphism $\varphi$ form $A$ into $M(B)$ is
nondegenerate (i.e. $B$ is spanned by vectors $\varphi(a)b $ and
spanned by vectors $b\varphi(a)$)then it has a unique extension to
a homomorphism \nolinebreak $M(A)\longrightarrow$\nolinebreak
$M(B)$.

\begin{itemize}
    \item A Hopf algebra over $\complex$ (the field of complex numbers) is a tuple $(A,\Delta,\varepsilon,S)$ where $A$ is an associative unital algebra over $\complex$, $\Delta:A\longrightarrow A\otimes A$,   $\varepsilon:A\longrightarrow \complex$ are algebra homomorphisms, and  $S:A\longrightarrow A$ is an algebra antihomomorphism such that $\Delta$ is coassociative,
$$
    (I\otimes \varepsilon)\Delta=(\varepsilon \otimes I)\Delta=I ,
$$
and
$$
m(I\otimes S)\Delta=m(S\otimes I)\Delta=\varepsilon 1_{A}
$$
where $m$ is the multiplication in $A$. For more details see \cite{1,2}.
\item A $\pi-$coalgebra over $\complex$ is a family $A=\{A_{\alpha}\}_{\alpha \in \pi}$($\pi$ is a discrete group) of $\complex$ spaces  endowed with a family $\Delta=\{\Delta_{\alpha,\beta}:A_{\alpha\beta}\longrightarrow A_{\alpha}\otimes A_{\beta}\}$ of $\complex-$linear maps (the comultiplication) and a linear map $\varepsilon:A_{1}\longrightarrow \complex$ (the counit) such that
\begin{enumerate}
    \item For any $\alpha,\beta,\gamma \in \pi\;; \Delta $ is coassociative in the sense that
    $$(\Delta_{\alpha,\beta}\otimes I)\delta_{\alpha\beta,\gamma}=(I\otimes \Delta_{\beta,\gamma})\Delta_{\alpha,\beta\gamma},$$
    \item For all $\alpha \in \pi$
    $$(I\otimes \varepsilon)\Delta_{\alpha,1}=(\varepsilon \otimes I)\Delta_{1,\alpha}=I$$
where 1 is the identity in in the group and $I$ is the identity map .
\end{enumerate}
\item A Hopf $\pi-$coalgebra  over $\complex$ is a $\pi-$coalgebra $H=(\{H_{\alpha}\}_{\alpha\in\pi},\Delta,\varepsilon)$ endowed with a family $S=\{S_{\alpha}:H_{\alpha}\longrightarrow H_{\alpha^{-1}}\}_{\alpha \in \pi}$ of $\complex-$linear maps (the antipode)such that
\begin{enumerate}
    \item Each $H_{\alpha}$ is an algebra with multiplication $m_{\alpha}$ and unit element $1_{\alpha} \in H_{\alpha}$,
    \item $\varepsilon:H_{1}\longrightarrow \complex$ and $\Delta_{\alpha,\beta}:H_{\alpha\beta}\longrightarrow H_{\alpha}\otimes H_{\beta}$ (for all $\alpha,\beta \in \pi $) are algebra homomorphisms ,
    \item For any $\alpha \in \pi$,
    $$m_{\alpha}(S_{\alpha^{-1}}\otimes I)\Delta_{\alpha^{-1},\alpha}=\varepsilon 1_{\alpha}=m_{\alpha}(I\otimes S_{\alpha^{-1}})\Delta_{\alpha,\alpha^{-1}}.$$
\end{enumerate}
    \item A multiplier Hopf algebra is a pair $(A,\Delta)$ where $A$ is an associative algebra over $\complex$, with or without identity and $\Delta :A\longrightarrow M(A\otimes A)$ is an algebra homomorphism, $M(A\otimes A)$ is the multiplier algebra of $A\otimes A$, such that
\begin{enumerate}
    \item  For all $a,b\in A $\\
    $$\Delta(a)(1\otimes b)\in A\otimes A \;\;\;\;\;\textrm{and} \;\;\;\;\; (a\otimes1)\Delta(b)\in A\otimes A,$$
    \item  $\Delta $ is coassociative in the sense that
    $$
(a\otimes 1\otimes 1)(\Delta \otimes I)(\Delta (b)(1\otimes c))=(I\otimes \Delta)((a\otimes 1)\Delta (b))(1\otimes 1\otimes c)$$
for all $a,b,c,\in A,$
    \item The linear maps
\begin{eqnarray}
&T^{1}:A\otimes A\longrightarrow A\otimes A&, \nonumber \\
&T^{2}:A\otimes A\longrightarrow A\otimes A& \nonumber
\end{eqnarray}
defined by
\begin{eqnarray*}
&T^{1}(a \otimes b)=\Delta(a)(1\otimes b)&,
\nonumber \\
&T^{2}(a \otimes b)=(a\otimes 1)\Delta(b)&
\nonumber
\end{eqnarray*}
are bijective.
\end{enumerate}
\end{itemize}

\section{Multiplier Hopf group coalgebra.}
Assume that $A=\{A_{\alpha}\}_{\alpha\in \pi}$ (where $\pi$ is a discrete group ) is a family of associative algebras over $\complex$ with nondegenerate products, $ A_{\alpha}$ may or may not have a unit.
\begin{defn}
 A comultiplication on $A$ is a family of
homomorphisms
$$
\Delta =\{\Delta_{\alpha,\beta}:A_{{\alpha\beta}}\longrightarrow
M(A_{\alpha} \otimes A_{\beta})\}_{\alpha,\beta \in \pi}
$$
such that for any $\alpha,\beta,\gamma \in \pi$
\begin{enumerate}
\item  $\Delta _{\alpha ,\beta }(a)(1\otimes b)\in A_{\alpha }\otimes
A_{\beta }\;\;\;\forall a\in A_{\alpha\beta},b\in A_{\beta}$\hspace{0.1in} and \hspace{0.1in}
$(a\otimes 1)\Delta _{\alpha ,\beta }(b)\in A_{\alpha
}\otimes A_{\beta }\;\;\;\forall a\in A_{\alpha},b\in A_{\alpha\beta}$,
\item  $\Delta $ is coassociative in the sense that
$$
(a\otimes 1\otimes 1)(\Delta _{\alpha ,\beta }\otimes I)(\Delta _{\alpha
\beta ,\gamma }(b)(1\otimes c))=(I\otimes \Delta _{\beta ,\gamma
})((a\otimes 1)\Delta _{\alpha ,\beta \gamma }(b))(1\otimes 1\otimes c)
$$
 for all $a\in A_{\alpha},\;b\in A_{\alpha\beta\gamma},$ and $c\in A_{\gamma}$.
\end{enumerate}
\end{defn}

Condition (1) makes sense because $A_{\alpha}\otimes A_{\beta}\subseteq M(A_{\alpha}) \otimes M(A_{\beta})\subseteq M(A_{\alpha}\otimes A_{\beta})$. Also, if $A_{\alpha}$ is unital algebra $\forall \alpha\in \pi$, the above definition gives rise to the definition of comultiplication introduced by Tureav\cite{7}.
%============================================================================================================
\begin{defn}
Let $A=\{A_{\alpha}\}_{\alpha\in \pi}$ be a family of
algebras with a nondegenerate products over $\complex$, $A_{\alpha}$ may or may not
have a unit and let $\Delta$ be a comultiplication on $A$ .We call $(A,\Delta)$ a
multiplier Hopf $\pi-$coalgebra if the linear maps
\begin{eqnarray}
&T_{\alpha,\beta}^{1}:A_{\alpha\beta}\otimes A_{\beta}\longrightarrow
A_{\alpha}\otimes A_{\beta}& ,  \nonumber \\
&T_{\alpha,\beta}^{2}:A_{\alpha}\otimes A_{\alpha \beta}\longrightarrow
A_{\alpha}\otimes A_{\beta}&  \nonumber
\end{eqnarray}
defined by
\begin{eqnarray}
&T_{\alpha,\beta}^{1}(a \otimes b)=\Delta_{\alpha,\beta}(a)(1\otimes b)&,
\nonumber \\
&T_{\alpha,\beta}^{2}(a \otimes b)=(a\otimes 1)\Delta_{\alpha,\beta}(b)&
\nonumber
\end{eqnarray}
are bijective for all $\alpha,\beta \in \; \pi$ .
\end{defn}
%==================================================================================================================

\begin{rem}
In fact the condition that $T_{\alpha ,\beta }^{1},T_{\alpha ,\beta }^{2}$
are bijective implies that $\Delta $ is nondegenerate homomorphism, and hence the
homomorphisms $(\Delta _{\alpha ,\beta }\otimes I)$ and $(I\otimes \Delta
_{\beta ,\gamma })$ have  unique extensions to $M(A_{\alpha \beta }\otimes
A_{\gamma })$ and $M(A_{\alpha }\otimes A_{\beta \gamma })$ respectively  so the coassociativity  will just mean
$$
(\Delta _{\alpha ,\beta }\otimes I)\Delta _{\alpha \beta ,\gamma }=(I\otimes
\Delta _{\beta ,\gamma })\Delta _{\alpha ,\beta \gamma }
$$
as maps from $A_{\alpha\beta\gamma}$ to $M(A_{\alpha}\otimes A_{\beta}\otimes A_{\gamma})$.\\
Also, we can see that $(A_{1},\Delta_{1,1})$ is a multiplier Hopf algebra.
\end{rem}

\begin{defn}
 A multiplier Hopf $\pi -$coalgebra $(A,\Delta )$ is called  regular multiplier Hopf $\pi -$coalgebra if the opposite comultiplication $\Delta ^{^{\prime }}=\sigma \Delta $ satisfies
\begin{enumerate}
\item  $\Delta _{\alpha ,\beta }^{^{\prime }}(a)(1\otimes b)\in A_{\beta }\otimes A_{\alpha },\,\,\,(a\otimes 1)\Delta _{\alpha ,\beta }^{^{\prime }}(b)\in A_{\beta}\otimes A_{\alpha}\,\,\,$ for all $\alpha ,\beta \in \pi$,
\item  $\Delta ^{^{\prime }}$ is opposite coassociative in the following sense
$$
(a\otimes 1\otimes 1)(\Delta _{\beta ,\alpha }^{^{\prime }}\otimes I)(\Delta
_{\gamma ,\beta \alpha }^{^{\prime }}(b)(1\otimes c))=(I\otimes \Delta
_{\gamma ,\beta }^{^{\prime }})((a\otimes 1)\Delta _{\gamma \beta ,\alpha
}^{^{\prime }}(b))(1\otimes 1\otimes c),
$$

\item  The linear maps
\begin{eqnarray*}
T_{\alpha ,\beta }^{^{\prime }1} &:&A_{\alpha \beta }\otimes A_{\alpha }\longrightarrow A_{\beta }\otimes A_{\alpha}, \\
T_{\alpha ,\beta }^{^{\prime }2} &:&A_{\beta }\otimes A_{\alpha \beta }\longrightarrow A_{\beta }\otimes A_{\alpha }
\end{eqnarray*}
defined by
\begin{eqnarray*}
T_{\alpha ,\beta }^{^{\prime }1}(a\otimes b) &=&\Delta _{\alpha ,\beta }^{^{\prime }}(a)(1\otimes b), \\
T_{\alpha ,\beta }^{^{\prime }2}(a\otimes b) &=&(a\otimes 1)\Delta _{\alpha ,\beta }^{^{\prime }}(b)
\end{eqnarray*}
are bijective.
\end{enumerate}
\end{defn}
\medskip
\noindent\textbf{Example.} Let $\pi$ be any group, and  $A$ be
the algebra of complex finitely supported functions on $\pi$. In
this case $M(A)$ consists of all complex functions on $\pi$.
Moreover $A\otimes A$ can be naturally identified with finitely
supported functions on $\pi\otimes \pi$ so that $M(A\otimes A)$
consists of
all complex functions on $\pi\otimes \pi$.\\
If we define $A_\pi=\{A_\alpha\}_{\alpha\in\pi}$ where
$A_{\alpha}=A$\hspace{.2in}for all $\alpha\in \pi$ ,and
$\Delta_{\alpha,\beta}:A_{\alpha\beta}\longrightarrow
M(A_{\alpha}\otimes A_{\beta})$ where
$$
\Delta_{\alpha,\beta}(f)(s,t)=f(\beta^{-1}s\beta t),
$$
we will clearly get a family of homomorphisms
$\Delta=\{\Delta_{\alpha,\beta}\}_{\alpha,\beta\in\pi}$. If $f\in
A_{\alpha\beta}, g\in A_{\beta}$ then $(s,t)\longrightarrow
f(\beta^{-1}s\beta t)g(t)$ will have finite support, and if $g\in
A_{\alpha}, f\in A_{\alpha\beta}$ then $(s,t)\longrightarrow
g(s)f(s\beta^{-1}t\beta)$ will have finite support. By direct
calculations we can see that $\Delta$ is coassociative hence it is
a comultiplication. The maps
$T_{\alpha,\beta}^{1},T_{\alpha,\beta}^{2}$ are bijective with
inverses $R_{\alpha,\beta}^{1},R_{\alpha,\beta}^{2}$ defined by
\begin{eqnarray*}
R^{1}_{\alpha,\beta}(f)(s,t)&=& f(\beta st^{-1}\beta^{-1},t),\\
R^{1}_{\alpha,\beta}(f)(s,t)&=& f(s,\beta s^{-1}t\beta^{-1}).
\end{eqnarray*}

Let us now browse the two approaches for studying the structure of
multiplier Hopf group coalgebra.
%====================================================3333333333333333333333333=======================================
\section{The first approach.}
In this section we prove the existence of a suitable counit and antipode for the multiplier Hopf group coalgebra. Also, we prove that there exist a large subspace of the dual space that can be made into an algebra.
Finally, we prove that $A$ is a Hopf group coalgebra if and only if it is a unital multiplier Hopf group coalgebra, and  that $S^{-1}$ exists under some conditions.
\subsection{Construction of the couint.}
\begin{defn}
For every $\alpha \in \pi $  define a linear map
$$
E_{\alpha}:A_{1}\longrightarrow L(A_{\alpha})
$$
by
$$
E_{\alpha}(a)b=\,m_{\alpha}\left(T_{1,\alpha}^{1}\right)^{-1}(a\otimes b) \;\;\; \forall \;\;b \in A_{\alpha},
$$
where $L(A_{\alpha})$ is the left
multiplier algebra of $A_{\alpha}$, $m_{\alpha}:A_{\alpha} \otimes
A_{\alpha}\longrightarrow A_{\alpha}$ is the multiplication on $A_{\alpha}
\;$ and $T_{1,\alpha}^{1}(a \otimes b)=\Delta_{1,\alpha}(a)(1\otimes b)$.
\end{defn}
 Because $T_{1,\alpha }^{1}(x(1\otimes c))=T_{1,\alpha }^{1}(x)(1\otimes c)$
for all $x\in A_{\alpha }\otimes A_{\alpha },c\in A_{\alpha }$, the same result
holds for $\left( T_{1,\alpha }^{1}\right) ^{-1}$. And since  $m_{\alpha
}(x(1\otimes c))=m_{\alpha }(x)c$ we get that $E_{\alpha }(a)$ is a left
multiplier on $A_{\alpha }$ for all $a\in A_{1},\alpha \in \;\pi $.\\
Following A.Van Daele we will show that $E_{\alpha}(A_{1})\subset \complex 1_{\alpha}$ which will give  the map $\varepsilon$.
%=================================================================================================================
\begin{lem}\label{1co}
For all $a,b\in A_{\alpha},\;\;\alpha\in \pi$ we have
$$(I \otimes E_{\alpha})((b\otimes 1)\Delta_{\alpha,1}(a))=ba\otimes 1. $$
\end{lem}
\begin{proof}
Suppose that $a,b \in A_{\alpha}$, and let
$$
a\otimes b =\sum_{i=1}^{n} \Delta_{\alpha,\alpha}(a_{i})(1\otimes b_{i}).
$$
If we apply $(\Delta_{\alpha,1}\otimes I)$ and multiply by $(c\otimes1\otimes 1)$ on the left then
\begin{eqnarray*}
(c\otimes 1)\Delta_{\alpha,1}(a)\otimes b &=& \sum (c\otimes 1\otimes 1)(\Delta_{\alpha,1}\otimes I) (\Delta_{\alpha,\alpha}(a_{i})(1\otimes b_{i}))  \\
&=& \sum (I\otimes \Delta_{1,\alpha})((c\otimes 1)\Delta_{\alpha,\alpha}(a_{i}))(1\otimes 1\otimes b_{i}).
\end{eqnarray*}
Now, let $\varphi$ be any linear functional on $A_{\alpha}$. If we apply $
(\varphi \otimes I \otimes I)$ to the above equation, we get
\begin{eqnarray*}
(\varphi \otimes I)((c \otimes 1)\Delta_{\alpha,1}(a))\otimes b &=&(\varphi \otimes I \otimes I) (\sum (I \otimes \Delta_{1,\alpha})((c\otimes 1) \Delta_{\alpha,\alpha}(a_{i}))(1\otimes1\otimes b_{i}))   \\
&=& \sum \Delta_{1,\alpha}((\varphi \otimes 1)((c\otimes 1)\Delta_{\alpha,\alpha}(a_{i})))(1\otimes b_{i}) \\
&=& \sum T_{1,\alpha}^{1}((\varphi \otimes I)((c\otimes 1)\Delta_{\alpha,\alpha}(a_{i}))\otimes b_{i}).
\end{eqnarray*}
By the definition of $E_{\alpha}$ we have
$$
E_{\alpha}((\varphi \otimes I)((c \otimes 1)\Delta_{\alpha,1}(a)))\;b = \sum
(\varphi \otimes I)((c \otimes 1) \Delta_{\alpha,\alpha}(a_{i}))\;b_{i}
$$
then
\begin{eqnarray*}
(\varphi \otimes I)((I\otimes E_{\alpha})((c\otimes 1) \Delta_{\alpha,1}(a))(1\otimes b)) &=&  (\varphi \otimes I)((c\otimes 1)\sum \Delta_{\alpha,\alpha}(a_{i})(1\otimes b_{i})) \\
&=&(\varphi \otimes I)((c\otimes 1)(a\otimes b)).
\end{eqnarray*}
Since this holds for all $\varphi \in\;A_{\alpha}^{^{\prime}}$ we get
$$
(I\otimes E_{\alpha})((c\otimes 1) \Delta_{\alpha,1}(a))(1\otimes b)= (ca
\otimes 1)(1\otimes b).
$$
This gives the required formula in the multiplier algebra.

\end{proof}
%====================================================================================================
\begin{lem}
$E_{\alpha}(A_{1})\subset \complex 1_{\alpha}\;\;\;\; \forall \alpha\in\pi$
\end{lem}

\begin{proof}
By the bijectivity of $T_{\alpha ,1}^{2}$ defined by $T_{\alpha ,1}^{2}(p\otimes q)=(p\otimes 1)\Delta _{\alpha ,1}(q)$ we have
$$
(I\otimes E_{\alpha })(a\otimes b)\in \;A_{\alpha }\otimes 1\;\;\forall
\;\;a\in \;A_{1},\;\;b\in \,A_{\alpha }.
$$
Which gives that $E_{\alpha }(A_{1})\subseteq \complex 1_{\alpha }$ for all $\alpha\in \pi$.

\end{proof}
%======================================================================================================
\begin{defn}
Define for every $\alpha\in \pi$ a map
$$\varepsilon_{\alpha}:A_{1}\longrightarrow \complex$$
by
$$\varepsilon_{\alpha}(a)1_{\alpha}=E_{\alpha}(a) \hspace{.4in}\forall a\in A_{1}.$$
\end{defn}
%===============================================================================================================
Now we can rewrite the formula in lemma \ref{1co}  as
$$(I\otimes \varepsilon _{\alpha })((a\otimes 1)\Delta _{\alpha ,1}(b))=ab.$$
By the definition of $\varepsilon_{\alpha} $ we have
$$
(\varepsilon _{\alpha }\otimes I)(a\otimes b)=m_{\alpha }\left( T_{1,\alpha }^{1} \right)^{-1}(a\otimes b)
$$
and hence
$$
(\varepsilon _{\alpha }\otimes I)(\Delta _{1,\alpha }(a)(1\otimes b)=ab.
$$
%=================================================================================================================
\begin{lem}
$\varepsilon _{\alpha }$ is a homomorphism for all $\alpha\in\pi$
\end{lem}
\begin{proof}
Let $a\in A_{\alpha}$ and $b,c\in A_{1}$. By the bijectivity of $T_{\alpha,1}^{1}$ we have
\begin{eqnarray*}
a\otimes bc&=&\sum_{i} (1\otimes b)(d_{i}\otimes 1)\Delta_{\alpha,1}(c_{i})\\
&=& \sum_{i,j} (a_{ij}\otimes 1) \Delta_{\alpha,1}(b_{ij})\Delta_{\alpha,1}(c_{i})\\
&=& \sum_{i,j} (a_{ij}\otimes 1)\Delta_{\alpha,1}(b_{ij}c_{i})
\end{eqnarray*}
where $a_{ij}, b_{ij}, c_{i}, d_{i} \in A_{\alpha}$. Applying $(I\otimes \varepsilon_{\alpha})$ we get
$$
a\varepsilon_{\alpha}(bc)=\sum_{i,j}a_{ij} b_{ij} c_{i}.
$$
On the other hand we have
$$
a\varepsilon_{\alpha}(c)=\sum_{i}d_{i} c_{i}
$$
if
$$
a\otimes c=\sum_{i} (d_{i}\otimes 1)\Delta_{\alpha,1}(c_{i}).
$$
Then
\begin{eqnarray*}
a\otimes b\varepsilon_{\alpha}(c)&=& \sum_{i} (d_{i}\otimes b)(c_{i}\otimes 1)\\
&=&\sum_{i,j}(a_{ij}\otimes 1) \Delta_{\alpha,1}(b_{ij})(c_{i}\otimes 1).
\end{eqnarray*}
Applying $(I\otimes \varepsilon_{\alpha})$ we get
$$
a\varepsilon_{\alpha}(b)\varepsilon_{\alpha}(c)=\sum_{i,j}a_{ij} b_{ij} c_{i}.
$$
This means
$$
a\varepsilon_{\alpha}(bc)=a\varepsilon_{\alpha}(b)\varepsilon_{\alpha}(c).
$$
\end{proof}
%============================================================================================================
\begin{lem}
For all $\alpha,\beta \in \pi$ \hspace{.2in}$\varepsilon_{\alpha}=\varepsilon_{\beta}$.
\end{lem}

\begin{proof}
We can see that the map $$\varphi
:E_{\alpha}(A_{1})\longrightarrow E_{\beta}(A_{1})$$ defined by $$
\varphi(E_{\alpha}(a))=E_{\beta}(a)\;\;\; \forall a\in A_{1}$$
 is an algebra homomorphism. Then
$$
\varepsilon_{\beta}(a)I_{\beta}=E_{\beta}(a)=\varphi(E_{\alpha}(a))=
\varphi(\varepsilon_{\alpha}(a)I_{\alpha})=\varepsilon_{\alpha}(a)
\varphi(I_{\alpha})=\varepsilon_{\alpha}(a)I_{\beta}.
$$
Since this holds for all $a\in A_{1}$, then $\varepsilon_{\alpha}=\varepsilon_{\beta}$ for all $\alpha,\beta\in\pi$.
\end{proof}
%===============================================================================================================
Now, we can summarize the results in this section in the following theorem(generalize \cite{4} theorem \nolinebreak 3.6).
\begin{thm}  \label{th1}
 Let $(A,\Delta)$ be a multiplier Hopf $\pi-$coalgebra. Then there exist a homomorphism $\varepsilon:A_{1} \longrightarrow C$ such that
\begin{itemize}
    \item[] $ (I \otimes \varepsilon)((a\otimes 1)\Delta_{\alpha,1}(b))=ab,$
    \item []$ (\varepsilon \otimes I)(\Delta_{1,\alpha}(a)(1\otimes b))=ab $
    \end{itemize}
for all $ a,b \in A_{\alpha},\; \alpha\in \pi$.
\end{thm}
Back to our example when $f\in A_{1}, g\in A_{\alpha}$ and $t\in \pi$
$$
\varepsilon(f)g(t)= (m_{\alpha}(T_{1,\alpha})^{-1}(f\otimes g))(t)=(f\otimes g )(\alpha t t^{-1}\alpha^{-1},t)=f(e)g(t)
$$
and hence $\varepsilon(f)=f(e)$.
%==========================================4444444444444444444======================================================
\subsection{Construction of the antipode. }
\begin{defn}
For every $\alpha \in\;\pi$ define a map
$$
S_{\alpha}:A_{\alpha}\longrightarrow L(A_{\alpha^{-1}})
$$
by
$$
S_{\alpha}(a)b=(\varepsilon_{\alpha} \otimes
I)(T^{1}_{\alpha,\,\alpha^{-1}})^{-1}(a\otimes b) \;\;\; \forall \;\;b \in A_{\alpha^{-1}}
.$$
\end{defn}
As before $S_{\alpha }(a)$ is a left multiplier on $A_{\alpha
^{-1}}$ for all $a\in A_{\alpha }$.
%=================================================================================================================
\begin{lem} \label{1an}
 For all $a\in A_{1},\;\; b,c \in A_{\alpha},\;\; \alpha\in \pi$, we have
$$
(I\otimes S_{\alpha^{-1}})((c\otimes 1)\Delta_{\alpha,\alpha^{-1}}(a))(1\otimes b)=(c\otimes1)(T_{1, \alpha} ^{1})^{-1} (a\otimes b).
 $$
 \end{lem}

\begin{proof}
Assume that $a \in A_{1},\,\,b\in A_{\alpha}$. Write
$$
a\otimes b=\sum_{i=1}^{n} \Delta_{\alpha\alpha^{-1},\alpha}(a_{i})(1\otimes b_{i}).
$$
Applying $(\Delta_{\alpha,\alpha^{-1}}\otimes I)$, and multiplying by $(c\otimes  1\otimes 1)$ on the left, we obtain
\begin{eqnarray*}
(c\otimes 1)\Delta_{\alpha,\alpha^{-1}}(a)\otimes b &=& \sum (c\otimes 1\otimes1)(\Delta_{\alpha,\alpha^{-1}} \otimes I)(\Delta_{\alpha\alpha^{-1},\alpha}(a_{i})(1\otimes b_{i})) \\
&=&\sum (I\otimes \Delta_{\alpha^{-1},\alpha})((c\otimes 1)\Delta_{\alpha,\alpha^{-1}\alpha}(a_{i}))(1\otimes 1\otimes b_{i}).
\end{eqnarray*}
Let $\varphi\in A_{\alpha}^{^{\prime}}$, and apply $(\varphi \otimes I \otimes I)$ then
\begin{eqnarray*}
(\varphi\otimes I)((c\otimes 1) \Delta_{\alpha,\alpha^{-1}}(a))\otimes b &=& \sum \Delta_{\alpha^{-1},\alpha} ((\varphi\otimes I) ((c\otimes 1)\Delta_{\alpha,\alpha^{-1}\alpha}(a_{i})))(1\otimes b_{i}) \\
&=& \sum  T_{\alpha^{-1},\alpha}^{1}((\varphi\otimes I)((c\otimes 1)\Delta_{\alpha,\alpha^{-1}\alpha} (a_{i}))\otimes b_{i}).
\end{eqnarray*}
By the definition of $S_{\alpha^{-1}}$ we have
\begin{eqnarray*}
S_{\alpha^{-1}}((\varphi\otimes I)((c\otimes 1)\Delta_{\alpha,\alpha^{-1}}(a))b&=& \sum (\varepsilon \otimes
I)((\varphi\otimes I)((c\otimes 1)\Delta_{\alpha,\alpha^{-1}\alpha}(a_{i}))\otimes b_{i}) \\
&=&\sum (\varphi\otimes I)((I\otimes \varepsilon )((c\otimes 1)\Delta_{\alpha,\alpha^{-1}\alpha}(a_{i}))\otimes b_{i}) \\
&=& (\varphi\otimes I)(\sum ca_{i}\otimes b_{i}),
\end{eqnarray*}
and hence
$$
(\varphi\otimes I)((I\otimes S_{\alpha^{-1}})((c\otimes 1)\Delta_{\alpha,\alpha^{-1}}(a))(1\otimes b))=(\varphi\otimes I)((c\otimes 1) (T_{1,\alpha}^{1})^{-1}(a\otimes b)).
$$
Since this is true for all $\varphi\in A_{\alpha}^{^{\prime}}$, we get
$$
(I\otimes S_{\alpha^{-1}})((c\otimes 1) \Delta_{\alpha,\alpha^{-1}}(a))(1\otimes b)=(c\otimes1)(T_{1,\alpha}^{1}) ^{- 1}(a\otimes b).
$$
\end{proof}
%================================================================================================================
\begin{lem}\label{2an}
For all $a\in A_{1},\;\; b,c \in A_{\alpha},\;\; \alpha\in \pi$, we have
$$
m_{\alpha}(I\otimes S_{\alpha^{-1}})((c\otimes 1)\Delta_{\alpha,\alpha^{-1}}(a))(1\otimes b)=c\varepsilon(a) b.
$$
\end{lem}

\begin{proof}
If we apply $m_{\alpha}$ to the equation in lemma \ref{1an}, we get
\begin{eqnarray*}
m_{\alpha}((I\otimes S_{\alpha^{-1}})((c\otimes 1)\Delta_{\alpha,\alpha^{-1}}(a))(1\otimes b)) &=&m_{\alpha}((c\otimes1)(T_{1,\alpha}^{1})^{-1}(a\otimes b)) \\
&=&c\varepsilon(a)b.
\end{eqnarray*}
\end{proof}
%=================================================================================================================
\begin{lem}\label{3an}
$S_{\alpha}(a)$ is a right multiplier for all $a\in A_{\alpha},\;\; \alpha\in \pi$, and satisfies
$$
m_{\alpha }((c\otimes 1)(S_{\alpha ^{-1}}\otimes I)\Delta _{\alpha^{-1},\alpha }(a)(1\otimes b))=c\varepsilon (a)b
$$
for all $a\in A_{1},\;\; b,c \in A_{\alpha},\;\; \alpha\in \pi$.
\end{lem}
\begin{proof}
 For every $\alpha \in \pi $ define
$$
S_{\alpha }^{^{\prime }}:A_{\alpha }\longrightarrow R(A_{\alpha ^{-1}})  \;\;\; \forall \;\;b \in A_{\alpha^{-1}}
$$
by
$$
bS_{\alpha }^{^{\prime }}(a)=(I\otimes \varepsilon )(T_{\alpha ^{-1},\alpha
}^{2})^{-1}(b\otimes a)
$$
where $R(A_{\alpha })$ is the right multiplier algebra. Directly we can show that $S_{\alpha }^{^{\prime }}(a)$ is a right multiplier on $A_{\alpha ^{-1}}$ for all $a\in A_{\alpha }$.
\noindent Completely similar as in lemma \ref{1an}  we can prove that
$$
(c\otimes 1)(S_{\alpha ^{-1}}^{^{\prime }}\otimes I)(\Delta _{\alpha
^{-1},\alpha }(a)(1\otimes b))=(T_{\alpha ,1}^{2})^{-1}(c\otimes a)(1\otimes
b).
$$
If we apply $m_{\alpha }$ to the above equation, we get the  formula in the statement of the
lemma with $S^{^{\prime }}$ instead of $S$ since
$$
m_{\alpha }((T_{\alpha ,1}^{2})^{-1}(c\otimes a))=c\varepsilon (a).
$$
Now, we will show that $S_{\alpha }^{^{\prime }}=S_{\alpha }$ for all $\alpha
\in \pi $.\\
By definition of $S^{'}$ we have
$$
aS_{\alpha ^{-1}}^{^{\prime }}(b)=\sum a_{i}\varepsilon (b_{i})
$$
if
$$
a\otimes b=\sum (a_{i}\otimes 1)\Delta _{\alpha ,\alpha ^{-1}}(b_{i}).
$$
Applying $(I\otimes S_{\alpha ^{-1}})$, and multiplying by $1\otimes c$ we get
$$
a\otimes S_{\alpha ^{-1}}(b)c=\sum (I\otimes S_{\alpha ^{-1}})((a_{i}\otimes
1)\Delta _{\alpha ,\alpha ^{-1}}(b_{i}))(1\otimes c).
$$
Applying $m_{\alpha }$ we have
$$
aS_{\alpha ^{-1}}(b)c=\sum a_{i}\varepsilon (b_{i})c=aS_{\alpha
^{-1}}^{^{\prime }}(b)c,
$$
which shows that $S_{\alpha }^{^{\prime }}=S_{\alpha }$for all $\alpha \in
\pi $, thus $S_{\alpha}(a)$ is a multiplier on $A_{\alpha^{-1}}$.
\end{proof}
%=============================================================================================================
\begin{lem} For all $a,b \in A_{\alpha},\;\; \alpha\in \pi$
$$S_{\alpha}(ab)=S_{\alpha}(b)S_{\alpha}(a).$$
\end{lem}
\begin{proof}
 Let $a,b\in A_{\alpha^{-1}}$, $c\in A_{\alpha}$. Then
 \begin{eqnarray*}
 c\otimes a b &=&\sum_{i}(e_{i}\otimes a)\Delta_{\alpha,\alpha^{-1}}(b_{i})\\
 &=& \sum_{i,j}(c_{ij}\otimes 1)\Delta_{\alpha,\alpha^{-1}}(a_{i,j})\Delta_{\alpha,\alpha^{-1}}(b_{i}).
 \end{eqnarray*}
 Applying $(I\otimes S_{\alpha^{-1}})$, multiplying by $(1\otimes d)$, and applying $m_{\alpha}$ we get
 \begin{eqnarray*}
 c S_{\alpha^{-1}}(a b) d &=& \sum_{i,j} m_{\alpha}((I\otimes S_{\alpha^{-1}})((c_{ij}\otimes 1)\Delta_{\alpha,\alpha^{-1}}(a_{i,j}))(1\otimes d))\\
&=& \sum_{i,j} c_{ij}\varepsilon(a_{ij}b_{i})d.
\end{eqnarray*}
 On the other hand, we have
 $$ c S_{\alpha^{-1}}(b)=\sum_{i}e_{i}\varepsilon(b_{i})$$
 if
 $$
 c\otimes  b =\sum_{i}(e_{i}\otimes 1)\Delta_{\alpha,\alpha^{-1}}(b_{i}).
 $$
 Then
 \begin{eqnarray*}
c S_{\alpha^{-1}}(b)\otimes a &=& \sum_{i}e_{i}\otimes a \varepsilon(b_{i})\\
&=& \sum_{i,j}(c_{ij}\varepsilon(b_{i})\otimes 1)\Delta_{\alpha,\alpha^{-1}}(a_{i,j}).
 \end{eqnarray*}
  Applying $(I\otimes S_{\alpha^{-1}})$, multiplying by $(1\otimes d)$, and applying $m_{\alpha}$ we get
 $$ cS_{\alpha^{-1}}(b)S_{\alpha^{-1}}(a)d= \sum_{i,j} c_{ij}\varepsilon(a_{ij})\varepsilon(b_{i})d.
$$
 Thus
 $$S_{\alpha^{-1}}(a b)=S_{\alpha^{-1}}(b)S_{\alpha^{-1}}(a)
 $$
 for all $a ,b \in A_{\alpha^{-1}}$ and hence $S$ is an antihomomorphism.
\end{proof}
%============================================================================================================
The following theorem gives the main result of this section [generalize \cite{4} theorem 4.6].
\begin{thm} \label{th2}
If $(A,\Delta)$ is a multiplier Hopf $\pi-$coalgebra then there exist an antihomomorphism $S=\{S_{\alpha}:A_{\alpha}\longrightarrow M(A_{\alpha^{-1}})\}_{\alpha \in \pi}$ such that for all $a\in A_{1}\;\;\;b,c\in A_{\alpha}\;\;\;\alpha \in\,\pi$
\begin{itemize}
    \item[] $ m_{\alpha}((I\otimes S_{\alpha^{-1}})((c\otimes1)\Delta_{\alpha,\alpha^{-1}}(a))(1\otimes b))=c\varepsilon(a)b$,
\item[] $ m_{\alpha }((c\otimes 1)(S_{\alpha ^{-1}}\otimes I)\Delta _{\alpha^{-1},\alpha }(a)(1\otimes b))=c\varepsilon (a)b$.
\end{itemize}
\end{thm}
%===============================================================================================================
Back to our example let $f\in A_{\alpha}, g\in A_{\alpha^{-1}}$, and $t\in \pi$. Then
\begin{eqnarray*}
(S_{\alpha}(f)g)(t)&=& (\varepsilon\otimes I)(T^{1}_{\alpha,\alpha^{-1}})^{-1}(f\otimes g)(t)=(T^{1}_{\alpha,\alpha^{-1}})^{-1}(f\otimes g)(e,t)\\
&=&(f\otimes g)(\alpha^{-1} t^{-1}\alpha, t)= f(\alpha^{-1} t^{-1}\alpha)g(t)
\end{eqnarray*}
and hence $S_{\alpha}(f)(t)=f(\alpha^{-1} t^{-1}\alpha)$.
%===============================================================================================================
\begin{thm}
 $A$ is a Hopf $\pi-$coalgebra if and only if $A$ is a  multiplier Hopf $\pi-$coalgebra with unital component.
\end{thm}

\begin{proof}
 If  $(A,\Delta)$ is a Hopf $\pi-$coalgebra then  $A_{\alpha}$ has an identity $\forall \alpha \in \pi $ which means that
$ M(A_{\alpha})= A_{\alpha} \;\;\; \forall \alpha \in \pi $ and according to remark (1) we have $\Delta$ is a comultiplication.  \\
Now we will prove that the linear maps
\begin{eqnarray*}
&T_{\alpha,\beta}^{1}:A_{\alpha\beta}\otimes A_{\beta}\longrightarrow A_{\alpha}\otimes A_{\beta}&,  \\
&T_{\alpha,\beta}^{2}:A_{\alpha}\otimes A_{\alpha \beta}\longrightarrow A_{\alpha}\otimes A_{\beta}&
\end{eqnarray*}
defined by
\begin{eqnarray*}
&T_{\alpha,\beta}^{1}(a \otimes b)=\Delta_{\alpha,\beta}(a)(1\otimes b)&,  \\
&T_{\alpha,\beta}^{2}(a \otimes b)=(a\otimes 1)\Delta_{\alpha,\beta}(b)&
\end{eqnarray*}
are bijective for all $\alpha,\beta \in \pi\;.$\\
Consider the following linear maps
\begin{eqnarray*}
&R_{\alpha,\beta}^{1}:A_{\alpha}\otimes A_{\beta}\longrightarrow A_{\alpha \beta}\otimes A_{\beta}&,  \\
&R_{\alpha,\beta}^{2}:A_{\alpha}\otimes A_{\beta}\longrightarrow A_{\alpha}\otimes A_{\alpha \beta}&
\end{eqnarray*}
which defined by
\begin{eqnarray*}
&R_{\alpha,\beta}^{1}(a\otimes b)=((I\otimes S_{\beta^{-1}})\Delta_{\alpha\beta,\beta^{-1}}(a))(1\otimes b)&,   \\
&R_{\alpha,\beta}^{2}(a\otimes b)=(a\otimes 1)((S_{\alpha^{-1}}\otimes I)\Delta_{\alpha^{-1},\alpha\beta}(b))&.
\end{eqnarray*}
By using the properties of $S$ and $\Delta$ we can prove that $R^{1}_{\alpha,\beta}$ and  $R^{2}_{\alpha,\beta}$ are the inverses of $T^{1}_{\alpha,\beta}$ and $T^{2}_{\alpha,\beta}$ respectively. If we use Sweedler's notation we get
\begin{eqnarray*}
T^{1}_{\alpha,\beta}R^{1}_{\alpha,\beta}(a\otimes b)&=& \sum_{(a)}T^{1}_{\alpha,\beta}(a_{(1,\alpha\beta)} \otimes S_{\beta^{-1}}(a_{(2,\beta^{-1})})\,b)   \\
&=& \sum _{(a)}a_{(1,\alpha)}\otimes a_{(2,\beta)} S_{\beta^{-1}}(a_{(3,\beta^{-1})})\,b  \\
&=& \sum _{(a)}a_{(1,\alpha)}\otimes m_{\alpha}(I \otimes S_{\beta^{-1}})\Delta_{\beta,\beta^{-1}} (a_{(2,1)})(1\otimes b)  \\
&=& \sum _{(a)}(a_{(1,\alpha)}\otimes \varepsilon(a_{(2,1)}) 1_{\beta})(1 \otimes b)=\,a\otimes b.
\end{eqnarray*}
Similarly, for $R_{\alpha,\beta}^{1}T_{\alpha,\beta}^{1},R_{\alpha,\beta}^{2}T_{\alpha,\beta}^{2}$ and $T_{\alpha,\beta}^{2}R_{\alpha,\beta}^{2}$. Therefore   $A$ is a multiplier Hopf $\pi-$coalgebra.\\
If $A$ is a unital multiplier Hopf $\pi-$coalgebra then $ M(A_{\alpha}\otimes A_{\beta})=A_{\alpha}\otimes A_{\beta} \; \forall \:\alpha,\beta \in \pi $ and by using remark (1) we have
$$
(\Delta _{\alpha ,\beta }\otimes I)\Delta _{\alpha \beta ,\gamma }=(I\otimes \Delta _{\beta ,\gamma })\Delta _{\alpha ,\beta \gamma }
$$
as a maps from $A_{\alpha\beta\gamma}$ into $A_{\alpha}\otimes A_{\beta}\otimes A_{\gamma}$ \\
Also, since $\varepsilon$ is a nondegenerate homomorphism then $\varepsilon \otimes I$ and $I\otimes \varepsilon$ have unique extensions to $M(A_{1}\otimes A_{\alpha})$ and $M(A_{\alpha}\otimes A_{1})$, respectively. So the counitary property just means
$$ (I\otimes \varepsilon)\Delta_{\alpha,1}= (\varepsilon \otimes 1)\Delta_{1,\alpha}=I\;\;\;\; \forall \alpha\in \pi. $$
Since $A$ is unital then $ M(A_{1}\otimes A_{\alpha})=A_{1}\otimes A_{\alpha}$ and $M(A_{\alpha}\otimes A_{1})=A_{\alpha}\otimes A_{1}$ and then we obtain  the definition of the counit introduced by Turaev.\\
If $A$ is unital then the antihomomorphism $S$ has the form
$$ S=\{S_{\alpha}:A_{\alpha}\longrightarrow A_{\alpha^{-1}}\}_{\alpha \in \pi}$$ such that
$$ m_{\alpha}(S_{\alpha^{-1}}\otimes I)\Delta_{\alpha^{-1},\alpha}=\varepsilon 1_{\alpha}=m_{\alpha}(I\otimes S_{\alpha^{-1}})\Delta_{\alpha,\alpha^{-1}}.$$
\end{proof}

%==================================================================================================================
\subsection{Regular multiplier Hopf group coalgebra.}
 In this subsection, we reformulate the concept of regular multiplier Hopf algebra  to obtain a bijective antipode in our structure using the technique  of proofs of A. Van Daele.
\begin{defn}
 A multiplier Hopf $\pi -$coalgebra $(A,\Delta )$ is called  regular multiplier Hopf $\pi -$coalgebra if the opposite comultiplication $\Delta ^{^{\prime }}=\sigma \Delta $ satisfies
\begin{enumerate}
\item  $\Delta _{\alpha ,\beta }^{^{\prime }}(a)(1\otimes b)\in A_{\beta }\otimes A_{\alpha },\,\,\,(a\otimes 1)\Delta _{\alpha ,\beta }^{^{\prime }}(b)\in A_{\beta}\otimes A_{\alpha}\,\,\,$ for all $\alpha ,\beta \in \pi$,
\item  $\Delta ^{^{\prime }}$ is opposite coassociative in the following sense
$$
(a\otimes 1\otimes 1)(\Delta _{\beta ,\alpha }^{^{\prime }}\otimes I)(\Delta
_{\gamma ,\beta \alpha }^{^{\prime }}(b)(1\otimes c))=(I\otimes \Delta
_{\gamma ,\beta }^{^{\prime }})((a\otimes 1)\Delta _{\gamma \beta ,\alpha
}^{^{\prime }}(b))(1\otimes 1\otimes c),
$$

\item  The linear maps
\begin{eqnarray*}
T_{\alpha ,\beta }^{^{\prime }1} &:&A_{\alpha \beta }\otimes A_{\alpha }\longrightarrow A_{\beta }\otimes A_{\alpha}, \\
T_{\alpha ,\beta }^{^{\prime }2} &:&A_{\beta }\otimes A_{\alpha \beta }\longrightarrow A_{\beta }\otimes A_{\alpha }
\end{eqnarray*}
defined by
\begin{eqnarray*}
T_{\alpha ,\beta }^{^{\prime }1}(a\otimes b) &=&\Delta _{\alpha ,\beta }^{^{\prime }}(a)(1\otimes b), \\
T_{\alpha ,\beta }^{^{\prime }2}(a\otimes b) &=&(a\otimes 1)\Delta _{\alpha ,\beta }^{^{\prime }}(b)
\end{eqnarray*}
are bijective.
\end{enumerate}
\end{defn}
%=================================================================================================================
\begin{lem} \label{r1}
 For all $\alpha \in \pi $ and $a,b \in
A_{\alpha}$ we have

\begin{itemize}
\item  $(\varepsilon \otimes I)((1\otimes a)\Delta _{1,\alpha }(b))=ab$,

\item  $(I\otimes \varepsilon )(\Delta _{\alpha ,1}(a)(b\otimes 1))=ab$.
\end{itemize}
\end{lem}

\begin{proof}
Define
$$
E_{\alpha }^{^{\prime }}:A_{1}\rightarrow L(A_{\alpha })
$$
by
$$
E_{\alpha }^{^{\prime }}(a)b=m_{\alpha }(T_{\alpha ,1}^{^{\prime
}1})^{-1}(a\otimes b).
$$
Analogous to lemma \ref{1co} and using the opposite coassociativity we can prove that
$$
(I\otimes E_{\alpha }^{^{\prime }})((c\otimes 1)\Delta _{1,\alpha
}^{^{\prime }}(a))=ca\otimes 1.
$$
Define
 $$\varepsilon_{\alpha} ^{^{\prime }}(a)1_{\alpha }=E_{\alpha }^{^{\prime }}(a).$$
 Then
$$
(I\otimes \varepsilon ^{^{\prime }})((c\otimes 1)\Delta _{1,\alpha
}^{^{\prime }}(a_{i}))=ca.
$$
By the definition of $\varepsilon_{\alpha}^{^{\prime }}$ we  have
$$
( \varepsilon _{\alpha }^{^{\prime }}\otimes I)(a\otimes b)=m_{\alpha
}(T_{\alpha ,1}^{^{\prime }1})^{-1}(a\otimes b)
$$
and hence
$$
( \varepsilon _{\alpha }^{^{\prime }}\otimes I)(\Delta _{\alpha,1 }^{^{\prime
}}(a)(b\otimes 1))=ab.
$$

If we apply $\sigma $ to the above  equations we get
\begin{eqnarray*}
(\varepsilon _{\alpha }^{^{\prime }}\otimes I)((1\otimes a)\Delta _{1,\alpha }(b)) &=&ab,   \\
(I\otimes \varepsilon _{\alpha }^{^{\prime }})(\Delta _{\alpha ,1}(a)(b\otimes 1)) &=&ab.
\end{eqnarray*}
which are the equations in the statement of the lemma with $\varepsilon_{\alpha}
^{^{\prime }}$ instead of $\varepsilon $. We will now prove that $%
\varepsilon_{\alpha} ^{^{\prime }}=\varepsilon $. \\
\noindent By theorem \ref{th1} we have
$$
c(\varepsilon _{\alpha }\otimes I)(\Delta _{1,\alpha }(a)(1\otimes
b))=cab=(\varepsilon _{\alpha }\otimes I)((1\otimes c)\Delta _{1,\alpha
}(a))b.
$$
Then
$$
(\varepsilon _{\alpha }\otimes I)((1\otimes c)\Delta _{1,\alpha
}(a))=(\varepsilon _{\alpha }^{^{\prime }}\otimes I)((1\otimes c)\Delta
_{1,\alpha }(a)).
$$
By the bijectivity of the map $\sigma T_{\alpha ,\beta }^{^{\prime }2}$ we have
$$
\varepsilon _{\alpha }=\varepsilon _{\alpha }^{^{\prime
}}\;\;\;\;\;\;\;\;\;\forall \;\;\alpha \in \pi.
$$
\end{proof}
%================================================================================================================
\begin{defn}
For every $\alpha \in \pi$ define a linear map
$$
S_{\alpha }^{^{\prime }}(a):A_{\alpha }\rightarrow M(A_{\alpha
^{-1}})
$$
by
$$
S_{\alpha }^{^{\prime }}(a)(b)=(\varepsilon \otimes I\;)(T_{\alpha
^{-1},\alpha }^{^{\prime }1})^{-1}(a\otimes b).
$$
We define $S^{^{\prime }}=\{S_{\alpha}^{^{\prime }}\}_{\alpha\in\pi}$.
\end{defn}
%============================================================================================================
We will show that $S^{^{\prime }}$ has the same properties of $S$ in the sense of the definition of the regular multiplier Hopf group coalgebra and  is the inverse of $S$.
\begin{lem} \label{r2}
 If $(A,\Delta)$ is a regular multiplier Hopf group coalgebra then $S$ is invertible with inverse $S^{'}$,  $S_{\alpha }(A_{\alpha })\subseteq A_{\alpha
^{-1}}$, and $S_{\alpha }^{^{\prime }}(A_{\alpha })\subseteq A_{\alpha ^{-1}}$ for\ all$\;\alpha \in \pi $.
\end{lem}

\begin{proof}
Analogous  to lemmas \ref{1an}, \ref{2an}, and  \ref{3an}; and using the opposite coassociativity we can proof
$$
m_{\alpha }((I\otimes S_{\alpha ^{-1}}^{^{\prime}})((c\otimes 1)\Delta
_{\alpha ^{-1},\alpha }(a))(1\otimes b))=c\varepsilon (a)b.
$$
Similarly, we can prove that $ S_{\alpha ^{-1}}^{^{\prime }}(a)$ is  a
right multiplier for all  $a\in A_{\alpha }\;,\alpha \in \pi \;$and satisfies
$$
m_{\alpha }((c\otimes 1)(S_{\alpha ^{-1}}^{^{\prime }}\otimes I))((\Delta
_{\alpha ,\alpha ^{-1}}^{^{\prime }}(a)(1\otimes b)))=c\varepsilon (a)b.
$$
Now by the definition of $S_{\alpha }$ we have
$$
S_{\alpha }(a)(b)=\sum \varepsilon _{\alpha }(a_{i})b_{i}
$$
if
$$
a\otimes b=\sum \Delta _{\alpha ,\alpha ^{-1}}(a_{i})(1\otimes b_{i}).
$$
Then
\begin{eqnarray*}
b\otimes a &=&\sum \Delta _{\alpha ^{-1},\alpha }^{^{\prime
}}(a_{i})(b_{i}\otimes 1), \\
b\otimes ac &=&\sum \Delta _{\alpha ^{-1},\alpha }^{^{\prime
}}(a_{i})(b_{i}\otimes c).
\end{eqnarray*}
Applying ($S_{\alpha ^{-1}}\otimes I)\;$and multiplying by $(d\otimes 1)$
\begin{eqnarray*}
dS_{\alpha ^{-1}}(b)\otimes ac &=&\sum (d\otimes 1)(S_{\alpha
^{-1}}^{^{\prime }}\otimes I)(\Delta _{\alpha ,\alpha ^{-1}}^{^{\prime
}}(a_{i})(b_{i}\otimes c)) \\
&=&\sum (d\otimes 1)(S_{\alpha ^{-1}}^{^{\prime }}(b_{i})\otimes
1)(S_{\alpha ^{-1}}^{^{\prime }}\otimes I)(\Delta _{\alpha ,\alpha
^{-1}}^{^{\prime }}(a_{i})(1\otimes c)).
\end{eqnarray*}
Applying $m_{\alpha }\;$
\begin{eqnarray*}
dS_{\alpha ^{-1}}(b)ac &=&\sum dS_{\alpha ^{-1}}^{^{\prime
}}(b_{i})\varepsilon (a_{i})c \\
&=&\sum dS_{\alpha ^{-1}}^{^{\prime }}(\varepsilon (a_{i})b_{i})c \\
&=&dS_{\alpha ^{-1}}^{^{\prime }}(S_{\alpha }(a)b)c.
\end{eqnarray*}
Which shows that
$$
S_{\alpha ^{-1}}^{^{\prime }}(b)a=S_{\alpha ^{-1}}^{^{\prime }}(S_{\alpha
}(a)b).
$$
Since the elements of the form $S_{\alpha }(a)b\;$generates $A_{\alpha
^{-1}}$ then the above formula implies that $S_{\alpha ^{-1}}^{^{\prime
}}(A_{\alpha ^{-1}})\subseteq A_{\alpha }\;$for all $\alpha \in \pi$.
 Similarly, $S_{\alpha ^{-1}}(A_{\alpha ^{-1}})\subseteq A_{\alpha }$. The
above formula gives that
$$
S_{\alpha ^{-1}}^{^{\prime }}(b)a=S_{\alpha ^{-1}}^{^{\prime }}(b)S_{\alpha
^{-1}}^{^{\prime }}(S_{\alpha }(a)).
$$
If we multiply by $c\;$ on the left and use that elements of the form $
cS_{\alpha ^{-1}}^{^{\prime }}(b)$ span $A$, we get $a=S_{\alpha
^{-1}}^{^{\prime }}(S_{\alpha }(a))$. Similarly, we can show that $b=S_{\alpha }(S_{\alpha
^{-1}}^{^{\prime }}(b))$.
\end{proof}
%===============================================================================================================
\begin{lem} \label{a}
\begin{enumerate}
\item  Let $a\in A_{\alpha \beta }\;,b\in A_{\beta }\;,a_{i}\in A_{\alpha
}$, and $b_{i}\in A_{\beta ^{-1}}$. Then
$$
a\otimes S_{\beta }(b)=\sum \Delta _{\alpha \beta ,\beta
^{-1}}(a_{i})(1\otimes b_{i})\Leftrightarrow \sum a_{i}\otimes S_{\beta
}^{-1}(b_{i})=(1\otimes b)\Delta _{\alpha ,\beta }(a).
$$
\item  Let $a\in A_{\alpha \beta }\;,b\in A_{\alpha }\;,a_{i}\in A_{\beta },
$ and $\;b_{i}\in A_{\alpha ^{-1}}$. Then
$$
S_{\alpha ^{-1}}(b)\otimes a=\sum (b_{i}\otimes 1)\Delta _{\alpha
^{-1},\alpha \beta }(a_{i})\Leftrightarrow \Delta _{\alpha ,\beta
}(a)(b\otimes 1)=\sum S_{\alpha }^{-1}(b_{i})\otimes a_{i}.
$$
\end{enumerate}
\end{lem}

\begin{proof}
Assume that
$$
a\otimes S_{\beta }(b)=\sum \Delta _{\alpha \beta ,\beta
^{-1}}(a_{i})(1\otimes b_{i}).
$$
Applying ($\Delta _{\alpha ,\beta }\otimes I)\;$and multiplying by$\;(c\otimes
1\otimes 1)$ on the left we get
\begin{eqnarray*}
(c\otimes 1)\Delta _{\alpha ,\beta }(a)\otimes S_{\beta }(b) &=&\sum
(c\otimes 1\otimes 1)(\Delta _{\alpha ,\beta }\otimes I)(\Delta _{\alpha
\beta ,\beta ^{-1}}(a_{i})(1\otimes b_{i})) \\
&=&\sum (I\otimes \Delta _{\beta ,\beta ^{-1}})((c\otimes 1)\Delta _{\alpha
,1}(a_{i}))(1\otimes 1\otimes b_{i}).
\end{eqnarray*}
Let $\varphi \in A_{\alpha }^{^{\prime}}$. Apply $(\varphi \otimes I\otimes
I)$, we get
\begin{eqnarray*}
(\varphi \otimes I)((c\otimes 1)\Delta _{\alpha ,\beta }(a))\otimes S_{\beta
}(b)) &=&\sum \Delta _{\beta ,\beta ^{-1}}((\varphi \otimes I)((c\otimes
1)\Delta _{\alpha ,1}(a_{i})))(1\otimes b_{i}) \\
&=&\sum T_{\beta ,\beta ^{-1}}^{^{\prime}1}((\varphi \otimes I)((c\otimes
1)\Delta _{\alpha ,1}(a_{i}))\otimes b_{i}).
\end{eqnarray*}
By the definition of $S_{\beta }$ we have
\begin{eqnarray*}
S_{\beta }((\varphi \otimes I)((c\otimes 1)\Delta _{\alpha ,\beta }(a))) S_{\beta }(b)) &=&\sum \varepsilon ((\varphi \otimes I)((c\otimes 1)\Delta
_{\alpha ,1}(a_{i})))\otimes b_{i} \\
&=&\sum (\varphi \otimes I)((I\otimes \varepsilon )((c\otimes 1)\Delta
_{\alpha ,1}(a_{i}))\otimes b_{i}).
\end{eqnarray*}
Then
$$
(\varphi \otimes I)((I\otimes S_{\beta })((c\otimes 1)\Delta _{\alpha ,\beta
}(a))(1\otimes S_{\beta }(b))) =\sum (\varphi \otimes I)(ca_{i}\otimes b_{i}).
$$
Since this holds for all $\varphi \in A_{\alpha}^{\prime}$ we get
$$
(I\otimes S_{\beta })(c\otimes b)\Delta _{\alpha ,\beta }(a) =\sum
(ca_{i}\otimes b_{i}).
$$
By cancelling $c\;$we get
$$
(I\otimes S_{\beta })(1\otimes b)\Delta _{\alpha ,\beta }(a) =\sum
(a_{i}\otimes b_{i}). $$
Applying $I\otimes S_{\beta}^{-1}$ we have
$$
(1\otimes b)\Delta _{\alpha ,\beta }(a) =\sum a_{i}\otimes S_{\beta
}^{-1}(b_{i}).
$$
Similarly, we can prove part  (2) .
\end{proof}
%===============================================================================================================
\begin{lem} \label{r3}
Let $a\in A_{\alpha ^{-1}}\;,b\in A_{\beta
^{-1}\alpha ^{-1}}\;,a_{i}\in A_{\alpha ^{-1}}$ and $b_{i}\in A_{\beta }$. Then the following statements are equivalent :
\begin{enumerate}
\item  $\Delta _{\beta ,\beta ^{-1}\alpha ^{-1}}(a)(1\otimes b)=\sum \Delta
_{\beta ,\beta ^{-1}\alpha ^{-1}}(a_{i})(b_{i}\otimes 1)$,
\item  $a\otimes S_{\alpha \beta }^{-1}(b)=\sum (a_{i}\otimes 1)\Delta
_{\alpha ^{-1},\alpha \beta }(b_{i})$.
\end{enumerate}
\end{lem}

\begin{proof}
Let
$$
\Delta _{\beta ,\beta ^{-1}\alpha ^{-1}}(a)(b\otimes 1)=\sum \Delta _{\beta
,\beta ^{-1}\alpha ^{-1}}(a_{i})(b_{i}\otimes 1).
$$
For $c\in $ $A_{\alpha \beta }$ write
$$
b_{i}\otimes S_{\alpha \beta }(c)=\sum _{k}\Delta _{\beta ,\beta ^{-1}\alpha
^{-1}}(p_{ik})(1\otimes q_{ik})
$$
for all $i$. Then
$$
\Delta _{\beta ,\beta ^{-1}\alpha ^{-1}}(a)(1\otimes bS_{\alpha \beta
}(c))=\sum_{i,k}\Delta _{\beta ,\beta ^{-1}\alpha
^{-1}}(a_{i}p_{ik})(1\otimes q_{ik})
$$
and by the  bijectivity of the map $T_{\beta ,\beta ^{-1}\alpha ^{-1}}^{1}$ we have
$$
a\otimes bS_{\alpha \beta }(c)=\sum_{i,k} a_{i}p_{ik}\otimes q_{ik}.
$$
On the other hand, by Lemma \ref{r2} we have
$$
(1\otimes c)\Delta _{\alpha ^{-1},\alpha \beta }(b)=\sum_{i,k} p_{ik}\otimes
S_{\alpha \beta }^{-1}(q_{ik}).
$$
Therefore
\begin{eqnarray*}
a\otimes cS_{\alpha \beta }^{-1}(b) &=&\sum_{i,k}a_{i}p_{ik}\otimes
S_{\alpha \beta }^{-1}(q_{ik}) \\
&=&\sum_{i}(a_{i}\otimes c)\Delta _{\alpha ^{-1},\alpha \beta }(b_{i}).
\end{eqnarray*}
By cancelling $c$ we get the required formula.
\end{proof}
%==============================================================================================================
Analogous to lemma \ref{a} we can prove the following lemma.
%===============================================================================================================
\begin{lem}\label{r4}
 \begin{enumerate}
\item  Let $a\in A_{\beta ^{-1}\alpha ^{-1}}\;,b\in A_{\beta ^{-1}}\;,a_{i}\in A_{\alpha ^{-1}}$, and $\;b_{i}\in A_{\beta }$. Then
$$
a\otimes S_{\beta ^{-1}}^{^{\prime }}(b)=\sum \Delta _{\beta ,\beta ^{-1}\alpha ^{-1}}^{^{\prime }}(a_{i})(1\otimes b_{i})\Leftrightarrow (1\otimes b)\Delta _{\beta ^{-1},\alpha ^{-1}}^{^{\prime }}(a)=\sum
a_{i}\otimes S_{\beta ^{-1}}^{^{\prime }-1}(b_{i}).
$$

\item  Let $b\in A_{\alpha ^{-1}},b_{i}\in A_{\alpha },a\in A_{\alpha ^{-1}\beta ^{-1}},$ and $a_{i}\in A_{\beta ^{-1}}$. Then
$$
S_{\alpha ^{-1}}^{^{\prime }}(b)\otimes a=\sum (b_{i}\otimes 1)\Delta _{\beta ^{-1}\alpha ^{-1},\alpha }^{^{\prime }}(a_{i})\;\;\Leftrightarrow \;\;\Delta _{\beta ^{-1},\alpha ^{-1}}^{^{\prime }}(a)(b\otimes 1)=\sum
S_{\alpha ^{-1}}^{^{\prime }-1}(b_{i})\otimes a_{i}.
$$
\end{enumerate}
\end{lem}
%=================================================================================================================
\begin{lem} \label{r5}
\begin{enumerate}
\item  Let $a,a_{i}\in A_{\alpha ^{-1}}\;,b\in A_{\beta }$, and $b_{i}\in
A_{\beta ^{-1}\alpha ^{-1}}$. Then
$$
\Delta _{\beta ,\beta ^{-1}\alpha ^{-1}}^{^{\prime }}(a)(1\otimes b)=\sum
\Delta _{\beta ,\beta ^{-1}\alpha ^{-1}}^{^{\prime }}(a_{i})(b_{i}\otimes
1)\Leftrightarrow a\otimes S_{\beta ^{-1}}^{^{\prime }-1}(b)=\sum
(1\otimes a_{i})\Delta _{\beta ^{-1},\alpha ^{-1}}^{^{\prime }}(b_{i}).
$$
\item  Let $b\in A_{\beta ^{-1}\alpha ^{-1}}\;,a\in A_{\alpha
^{-1}}\;,a_{i}\in A_{\alpha ^{-1}}$, and $\;b_{i}\in A_{\beta }.$ Then
$$
\Delta _{\beta ,\beta ^{-1}\alpha ^{-1}}(a)(1\otimes b)=\sum \Delta _{\beta
,\beta ^{-1}\alpha ^{-1}}(a_{i})(b_{i}\otimes 1)\Leftrightarrow (1\otimes
a)\Delta _{\beta ^{-1},\alpha ^{-1}}(b)=\sum S_{\beta }(b_{i})\otimes a_{i}.
$$
\end{enumerate}
\end{lem}
%===========================================================================================================
The first part of the above lemma is a reformulation of lemma \ref{r3} for the opposite \nolinebreak comultiplication .
For part (2) if we apply $\sigma$ to the equivalence in part (1), we get
$$\Delta _{\beta ,\beta ^{-1}\alpha ^{-1}}(a)(b\otimes 1)=\sum
\Delta _{\beta ,\beta ^{-1}\alpha ^{-1}}(a_{i})(1\otimes
b_{i})\Leftrightarrow S_{\beta }(b)\otimes a =\sum
( a_{i}\otimes 1)\Delta _{\beta ^{-1},\alpha ^{-1}}(b_{i}).
$$
But this is up to a different summation gives the equivalence in part (2).

%==============================================================================================================
\begin{prop}
Let $b\in A_{\alpha },a\in A_{\alpha \beta }$. Then
$$
(1\otimes S_{\alpha }(b))\Delta _{\beta ^{-1},\alpha ^{-1}}(S_{\alpha \beta
}(a))=(S_{\beta }\otimes S_{\alpha })[\Delta _{\alpha ,\beta }^{'
}(a)(1\otimes b)].
$$
\end{prop}

\begin{proof}
Write
$$
\Delta _{\alpha ,\beta ^{-1}\alpha ^{-1}}(S_{\alpha }(b))(1\otimes S_{\alpha
\beta }(a))=\sum \Delta _{\beta ,\beta ^{-1}\alpha
^{-1}}(a_{i})(b_{i}\otimes 1).
$$
By lemma \ref{r5} we have
$$
S_{\alpha }(b)\otimes a=\sum (a_{i}\otimes 1)\Delta _{\alpha ^{-1},\alpha
\beta }(b_{i}),
$$
and
$$
(1\otimes S_{\alpha }(b))\Delta _{\beta ^{-1},\alpha ^{-1}}(S_{\alpha \beta
}(a))=\sum S_{\beta }(b_{i})\otimes a_{i}.
$$
By lemma \ref{r2} we have
$$
\Delta _{\alpha \beta }(a)(b\otimes 1)=\sum S_{\alpha }^{-1}(a_{i})\otimes b_{i}.
$$
If we combine all of this we get
\begin{eqnarray*}
(1\otimes S_{\alpha }(b))\Delta _{\beta ^{-1},\alpha ^{-1}}(S_{\alpha \beta
}(a)) &=&\sum (S_{\beta }\otimes S_{\alpha })(b_{i}\otimes S_{\alpha
}^{-1}(a_{i})) \\
&=&(S_{\beta }\otimes S_{\alpha })(\sigma (\Delta _{\alpha ,\beta
}(a)(b\otimes 1))) \\
&=&(S_{\beta }\otimes S_{\alpha })(\Delta _{\alpha ,\beta }^{'
}(a)(1\otimes b)).
\end{eqnarray*}
\end{proof}
%=============================================================================================================
\subsection{The Dual Algebra}
In general, if $\Delta(A)\subset A\otimes A $ then the dual space becomes an algebra, but this is not the case  in the multiplier Hopf $\pi-$coalgebra. So, in the following proposition we will use the technique of A. Van Daele in constructing the dual space.

\begin{prop}
If $(A,\Delta)$ is a multiplier Hopf $\pi-$coalgebra then the vector space $A^{\star}$ defined by
$$ A^{\star}=\oplus_{\alpha \in \pi } A^{\star}_{\alpha}, $$
where $A_{\alpha }^{*}$ is the vector space of linear functionals on $A_{\alpha }$
spanned by elements of the form $a\rightarrow f(bac)$ where $b,c\in
A_{\alpha }$ and $f\in A_{\alpha }^{\prime}$, can be made into associative algebra under  the product
$$
(fg)(a)=(f\otimes g)(\Delta _{\alpha ,\beta }(a))
$$
where $f\in A_{\alpha }^{*}$ and $g\in A_{\beta }^{*}.$
\end{prop}
\begin{proof}
 Let $f\in A^{*}_{\alpha}$, $g\in A^{*}_{\beta}$. Then there exist $b,c\in A_{\alpha}$, $d,e\in A_{\beta}$,  $f^{'}\in A^{'}_{\alpha}$, and  $g^{'}\in A^{'}_{\beta}$ such that  $f(a)=f^{\prime }(bac)$ for all $a\in A_{\alpha}$, and $g(a)=g^{'}(dae)$ for all $a\in A_{\beta}$. Assume that $a\in A_{\alpha,\beta}$ then
\begin{eqnarray*}
(fg)(a) &=&m_{c}(f\otimes g)(\Delta _{\alpha ,\beta }(a)) \\
&=&(f^{\prime}\otimes g^{\prime})((b\otimes d)\Delta _{\alpha
,\beta }(a)(c\otimes e)).
\end{eqnarray*}
Since $(b\otimes d)\Delta _{\alpha ,\beta }(a)(c\otimes e)\in $ $A_{\alpha
}\otimes A_{\beta }\;$then $(fg)(a)\in \complex\;$. It means that $fg$ is a linear
functional on $A_{\alpha \beta }.$ By linearity it can be extended to all $
f,g\in A^{*}.$\\
\indent By the bijectivity of $T_{\alpha ,\beta }^{1}\;$and$\;T_{\alpha ,\beta
}^{2}\;$we can write
\begin{eqnarray*}
b\otimes d &=&\sum_{i=1}^{n}(p_{i}\otimes 1)\Delta _{\alpha ,\beta }(q_{i}), \\
c\otimes e &=&\sum_{j=1}^{m} \Delta _{\alpha ,\beta }(r_{j})(1\otimes k_{j}).
\end{eqnarray*}
Then
$$
(fg)(a)=\sum_{i,j}(f^{\prime}\otimes g^{\prime})((p_{i}\otimes 1)\Delta _{\alpha
,\beta }(q_{i}ar_{j})(1\otimes k_{j})).
$$
Define
$$
h_{i,j}:A_{\alpha \beta }\rightarrow \complex
$$
by
$$
h_{i,j}(a)=(f^{\prime}\otimes g^{\prime})((p_{i}\otimes 1)\Delta _{\alpha
,\beta }(a)(1\otimes k_{j})).
$$
Then
$$
(fg)(a)=\sum_{i,j}h_{i,j}(qar).
$$

i.e. $fg\in A^{*}\;\;\; \forall f,g\in A^{\star}$
and thus $(A^{*},.)$ be an algebra.\\
Let $f_{1}\in A_{\alpha }^{*},f_{2}\in A_{\beta }^{*},$ and $f_{3}\in
A_{\gamma }^{*}\;$ and they defined by $f_{i}=f_{i}(b_{i}a_{i}c_{i})\; i=1,2,3.\;$where $a_{1},b_{1},c_{1}\in A_{\alpha },a_{2,}b_{2},c_{2}\in A_{\beta }\,
$and$\,a_{3},b_{3},c_{3}\in A_{\gamma }$. Then

\begin{eqnarray*}
((f_{1}f_{2})f_{3})(a) &=&(f_{1}f_{2}\otimes f_{3})((1\otimes b_{3})\Delta
_{\alpha \beta ,\gamma }(a)(1\otimes c_{3})) \\
&=&(f_{1}\otimes f_{2}\otimes f_{3})((b_{1}\otimes b_{2}\otimes 1)((1\otimes
b_{3})\Delta _{\alpha \beta ,\gamma }(a)(1\otimes c_{3}))(c_{1}\otimes
c_{2}\otimes 1) \\
&=&(f_{1}\otimes f_{2}\otimes f_{3})((1\otimes b_{2}\otimes
b_{3})((b_{1}\otimes 1\otimes 1)(\Delta _{\alpha ,\beta }\otimes I)(\Delta
_{\alpha \beta ,\gamma }(a)(1\otimes c_{3}))) \\
&&(c_{1}\otimes c_{2}\otimes 1)) \\
&=&(f_{1}\otimes f_{2}\otimes f_{3})((1\otimes b_{2}\otimes b_{3})((I\otimes
\Delta _{\beta ,\alpha })((b_{1}\otimes 1)(\Delta _{\alpha ,\beta \gamma
}(a))(1\otimes 1\otimes c_{3})) \\
&&(c_{1}\otimes c_{2}\otimes 1)) \\
&=&(f_{1}\otimes f_{2}\otimes f_{3})((1\otimes b_{2}\otimes b_{3})((I\otimes
\Delta _{\beta ,\alpha })((b_{1}\otimes 1)(\Delta _{\alpha ,\beta \gamma
}(a))(c_{1}\otimes 1)) \\
&&(1\otimes c_{2}\otimes c_{3})) \\
&=&(f_{1}(f_{2}f_{3}))(a).
\end{eqnarray*}
Then $A^{*}$ becomes an associative algebra.
\end{proof}
%=======================================================================================================
If $A_{\alpha}$ is a unital algebra $\forall \alpha\in\pi$, then we have the the same definition of the dual space introduced by Virelizier \cite{8}.
\begin{rem}
 If $f(a)=f^{\prime }(abc),f\in A_{\alpha }^{*}$, and if $\varepsilon$ is the counit of $A$ we have
\begin{eqnarray*}
\varepsilon f(a) &=&(\varepsilon \otimes f^{\prime })((1\otimes b)\Delta
_{1,\alpha }(a)(1\otimes c)) \\
&=&f^{\prime }(b(\varepsilon \otimes I)(\Delta (a)(1\otimes c))) \\
&=&f^{\prime }(bac) =f(a).
\end{eqnarray*}
Similarly, $f\varepsilon (a)=f(a)$.
So $\varepsilon $ corresponds to the identity in the multiplier algebra of $A^{*}.$
\end{rem}

\begin{defn}
 Let  $A$ be a  regular multiplier Hopf $\pi-$coalgebra and $S$ is the antipode of $A$. Then we define the adjoint antipode to be
$$
S^{\star }:A^{*}\rightarrow A^{*}
$$
by
$$
(S^{\star}f)(a)=f(S_{\alpha ^{-1}}(a))
$$
for all $f\in A_{^{\alpha }}^{*},a\in A_{^{\alpha ^{-1}}}$.
\end{defn}
%=========================================================================================================
\begin{lem}
$S^{\star }f\in A^{*}$ if $f\in A^{*}$.
\end{lem}
\begin{proof}
Let $f(a)=f^{' }(bac)$ where $f\in A_{\alpha }^{*},a,b,c\in A_{^{\alpha ^{-1}}}$, and choose $d\in A_{1}$ such that $\varepsilon (d)=1$. By using theorem \ref{th2} we get
\begin{eqnarray*}
(S^{^{\star }}f)(a) &=&f(S_{\alpha ^{-1}}(a))=f^{\prime }(bS_{\alpha
^{-1}}(a)c) \\
&=&f^{\prime }(b\varepsilon _{\alpha }(d)S_{\alpha ^{-1}}(a)c) \\
&=&f^{\prime }(m_{\alpha }((I\otimes S_{\alpha ^{-1}})((b\otimes 1)\Delta
_{\alpha ,\alpha ^{-1}}(d))(1\otimes S_{\alpha ^{-1}}(a)c))).
\end{eqnarray*}
By the  bijectivity of $T_{\alpha ,\alpha ^{-1}}^{2}\exists $ $p_{i}\otimes
q_{i}\in A_{\alpha }\otimes A_{^{\alpha ^{-1}}}$such that $(b\otimes 1)\Delta
_{\alpha ,\alpha ^{-1}}(d)=\sum_{i=1}^{n} p_{i}\otimes q_{i}$. Then
\begin{eqnarray*}
(S^{^{\star }}f)(a) &=&\sum f^{\prime }(m_{\alpha }(p_{i}\otimes S_{\alpha ^{-1}}(aq_{i})c) \\
&=&\sum f^{\prime }(p_{i}S_{\alpha ^{-1}}(aq_{i})c).
\end{eqnarray*}
By a similar argument we can prove that
$$
(S^{\star}f)(a)= \sum f^{^{\prime }}(b_{i}S_{\alpha ^{-1}}(r_{j}aq_{i})k_{j}).
$$
\end{proof}
%=============================================================================================================
\begin{lem}
If ($A,\Delta $) is a regular multiplier Hopf $\pi $- coalgebra then $S^{\star }$ is an antihomomorphism of $A^{*}$.
\end{lem}
\begin{proof}
Let $f\in A_{\beta ^{-1}}^{*},g\in A_{\alpha ^{-1}}^{*}$, and write $g(a_{1})=g^{\prime }(S_{\alpha }(b)a_{1})$ for some $b\in A_{\alpha }$.
Then for  $a\in A_{\alpha \beta }$ we have
\begin{eqnarray*}
S^{\star}(fg)(a) &=&fg(S_{\alpha \beta }(a))=(f\otimes g^{\prime })((1\otimes
S_{\alpha }(b))\Delta _{\beta ^{-1},\alpha ^{-1}}(S_{\alpha \beta }(a))) \\
&=&(f\otimes g^{\prime })((S_{\beta }\otimes S_{\alpha })(\Delta _{\alpha
\beta }^{\prime }(a)(1\otimes b))) \\
&=&(S^{\star}g^{\prime }\otimes S^{\star}f)(\Delta (a)(b\otimes 1)).
\end{eqnarray*}
But
$$
S^{\star}g^{\prime }(ab)=g^{\prime }(S_{\alpha }(b)S_{\alpha }(a))=g(S_{\alpha
}(a))=S^{\star}g(a),
$$
then
$$
S^{\star}(fg)=((S^{\star}g)(S^{\star}f))(a).
$$
\end{proof}
%=================================================================================================================
\begin{lem}
Let $A=\{A_p\}_{p\in G}$ be a multiplier Hopf $G-$coalgebra.  Then $H=\{p\in G :
A_p\neq 0\}$ is a subgroup of $G$.
\end{lem}
\begin{proof}
Let $H=\{p\in G : A_p\neq 0\}$ and  $p\in H$. Then from
the bijectivity of the map $T_{1,p}^{1}: A_p \otimes A_p\to
A_{1} \otimes A_p$ we get that $A_{1}\neq 0$. Also, if
$p,q\in H$ then the bijectivity of the map
$T_{p,q}^{1}: A_{p,q}\otimes A_q \to A_p\otimes A_q$
implies that $A_{p,q}\neq 0$. Finally, for all $p \in H$ the
bijectivity of the map $T_{p^{-1},p}^{1}:A_{1}\otimes A_p\to
A_{p^{-1}}\otimes A_p$ implies that $p^{-1}\in H$.
\end{proof}
%==================================================================================================================
\section{The second approach}
%==================================================================================================================
Now let $(A=\{A_p\}_{p\in G},\de=\{\de_{p,q}\}_{p,q\in G})$ be a \mhgc. Let $\ag$ be the direct sum of the algebras $A_p$. An element in $\ag$ is a function $a$ from $G$ to $\cup A_p$, with finite support and so that $a(p)\in A_p$ $\forall p\in G$. This algebra is an algebra with a nondegenerate product. For $a\in A_p$ we will use $a$ for the function in $\ag$ that have the value $a$ at $p$  and $0$ otherwise and use $a(p)$ when we mean the element  $a\in A_p$.
%=========================================================================================================
\begin{prop}
let $(A=\{A_p\}_{p\in G},\de=\{\de_{pq}\}_{p,q\in G})$ be a \mhgc. Then $\ag$ becomes a multiplier Hopf algebra under the comultiplication $\dg:\ag\to M(\ag\otimes\ag)$ given by
$$
\dg(a)(p,q)=\de_{p,q}(a(pq)).
$$
Moreover there exist a nondegenerate homomorphism $\ga: K(G) \to M(\ag)$ ,where $K(G)$ be the algebra of complex valued functions on $G$ with finite support,  such that $\ga(K(G))$ is in the center of  $M(\ag)$. And $\ga$ is compatible with the comultiplication.
\end{prop}
\begin{proof}
$I-$ $(A_{G},\Delta_{G})$ is a multiplier Hopf algebra. \\
$(1)$ Let $a,b \in \ag$. Then
$$
(\Delta_G(a)(1 \otimes b))(p,q)= \Delta_{p,q}(a(pq))(1 \otimes b(q)),
$$
but $b$ has a finite support, then $\Delta_G(a)(1 \otimes b) \in \ag \otimes \ag$. Similarly, $(a \otimes 1) \Delta_G(b) \in \ag \otimes \ag$. \\
$(2)$ Let $a,b,c \in \ag$. Then we have
\begin{eqnarray*}
(a \otimes 1 \otimes  1)(\dg \otimes I)(\dg(b)(1\otimes c))(p,q,r)&=&
(a(p) \otimes 1 \otimes 1)(\Delta_{p,q} \otimes I)[(\dg (b)(1 \otimes c))(pq,r)] \\
&=&(a(p) \otimes 1 \otimes 1) (\Delta_{p,q} \otimes I)[\Delta_{pq,r}(b(pqr))(1 \otimes c(r))]\\
&=&(I \otimes \Delta_{q,r})[(a(p) \otimes 1) \Delta_{p,qr}(b(pqr))](1 \otimes 1 \otimes c(r)) \\
&=&[(I \otimes \dg)[(a \otimes 1)\dg(b))](1 \otimes 1 \otimes c)](p,q,r).
\end{eqnarray*}
$(3)$ The map $T_{1}: A \otimes A \longrightarrow A \otimes A$ defined by
$$
T_{1}(a \otimes b)_{(p,q)}= \dg(a)(1 \otimes b)(pq)=\Delta_{p,q}(a(pq))(1 \otimes b (q))= T^{1}_{p,q}(a(pq) \otimes b(q))
$$
is a bijective with the inverse
$$
T^{-1}_{1}(a \otimes b)(p,q)=(T^{1}_{p,q})^{-1}(a(pq^{-1}) \otimes b(q)),
$$
since
\begin{eqnarray*}
((T^{-1}_{1} \circ T_{1})(a \otimes b))(p,q) &=& (T_{p,q}^{1})^{-1}(T_{1}(a \otimes b)(pq^{-1},q))\\
&=& (T^{1}_{p,q})^{-1}(T_{p,q}^{1}(a(pq^{-1}q) \otimes b(q))\\
&=& a(p) \otimes b(q).
\end{eqnarray*}
Similarly, $[(T_{1} \circ T_{1}^{-1}(a \otimes b)](p,q)= a(p) \otimes b(q)$.\\
Also, the map $T_{2}: A \otimes A \longrightarrow A \otimes A$ defined by
$$
T_{2}(a \otimes b)(p,q)=(a(p) \otimes 1)\Delta_{p,q}(b(pq))
$$
is bijective with the inverse
$$
T^{-1}_{2}(a \otimes b)(p,q)=(T^{2}_{p,q})^{-1}(a(p) \otimes b(p^{-1}q)).
$$
$II-$ There exist a nondegenerate homomorphism $\ga$, from $K(G)$ into the center of $M(\ag)$, compatible with the comultiplication.\\
$(1)$ Let $K(G)$ be the algebra of complex valued functions on the group $G$ with finite support.  Define
$$
\gamma: K(G) \longrightarrow M(A)
$$
by
$$
\gamma(f)(p)=f(p)1_{p}
$$
By the definition of $\ga$ we have that  $\gamma(f)$ in the center $M(A)\;\forall f \in K(G)$. Moreover, $\gamma(fg)(p)=(fg)(p)1_{p}=f(p)g(p)1_{P} = (\gamma(f)\gamma(g))(p)$ and $ \gamma(K(G))\ag=\ag\gamma(K(G))=\ag$.\\
$(2)$ Consider
\begin{eqnarray*}
( \gamma \otimes \gamma)\Delta_{\circ}(f)(p,q)&=&[(\gamma \otimes  \gamma)(\sum_{s}C_{s} \sum_{r}\delta_{r} \otimes  \delta_{r^{-1}s}(p,q))]\\
&=& [\sum_{s,r \in G}C_{s}\delta_{r}(p) \otimes  \delta_{r^{-1}s}(q) 1_{q}] \\
&=& f(pq)( 1_{p} \otimes 1_{q} ).
\end{eqnarray*}
On the other hand
\begin{eqnarray*}
\dg(\ga(f))(p,q)&=& \de_{p,q}(\ga(f)(p,q))\\
&=& \de_{p,q}(f(pq)1_{p,q})\\
&=& f(pq)( 1_{p} \otimes 1_{q} ).
\end{eqnarray*}
Thus $\ga$ is compatible with the comultiplication
\end{proof}
%===============================================================================================================
Let $(A,\Delta)$ be a multiplier Hopf algebra such that there
exist a nondegenerate homomorphism $\gamma: K(G)
\longrightarrow M(A)$ such that
$\gamma(K(G))$ is in the center of $M(A)$ and
$$
\Delta(\gamma(f))=( \gamma \otimes \gamma ) \Delta_{\circ}(f),
$$
where $\Delta_{\circ}(f)(p,q)=f(pq)$.\\
We denote this \mha with the nondegenerate homomorphism $\ga$ by $(\at,\dt)$ .\\

Let $A_{p}=\at \gamma (\delta_{p})$, $\delta _{p} \in K(G)$
where $\delta_{p}(q)= \delta^{q}_{p}= \{^{0, \;\;\;p \neq q}_{1,
\;\;\;p= q}$. Because of $A_p$ is a subalgebra of $\at$ and satisfies that if $a \gamma(\delta_{p}) b \gamma(\delta_{p})=0\; \forall \; a
\gamma(\delta_{p}),  \gamma(\delta_{p}) \neq 0 $ then $ab\gamma(\delta_{p})=0 \; \forall  \;a\in \at$ and hence $b\gamma(\delta_{p})=0$. Similarly,  if $a\gamma(\delta_{p}) b \gamma(\delta_{p}) =0 \; \forall
\;  b \gamma(\delta_{p}),  \gamma(\delta_{p}) \neq 0
$ then $a \gamma(\delta_{p})=0$.
 Then $A_{p}$ is an algebra with a nondegenerate
product. $\at$ can be identified with the direct of $A_{p}$'s using
$$
a \longrightarrow (a \gamma (\delta_{p}))_{p \in G}
$$
Since $\gamma $ is a nondegenerate homomorphism then for all $ a \in
\at$  there exist $f \in K(G) $ such that $a=a \gamma (f)$. Since $f \in
K(G)$ then there exist finite number of $\del's$ such that $f= \sum c_{p}
\delta_{p}$ and hence $a \in \oplus_{p}A\gamma(\delta_{p})=
\oplus_{p}A_{p}$. Then  $\forall a \in \at $ there exist a finite
numbers of $p$'s such that $a \in \oplus_{p}A_{p}$\\

Let $a,b \in \at$ and $\delta_{p}, \delta_{q} \in K(G)$. Then the map $\ttt_1$ satisfies
\begin{eqnarray*}
\ttt_{1}(a \gamma(\delta_{p}) \otimes b \gamma
(\delta_{q}))&=&\dt(a \gamma(\delta_{p}))(1 \otimes b \gamma(
\delta_{q}))\\
&=&\dt (a)( \gamma \otimes \gamma )
\Delta_{\circ}(\delta_{p})(1
\otimes b)(1 \otimes \gamma(\delta_{q})) \\
&=&\dt (a)(1 \otimes b)( \gamma \otimes \gamma
)T^{\circ}_{1}(\delta_{p} \otimes \delta_{q})\\
&=&\ttt_{1}(a \otimes b)( \gamma \otimes \gamma
)T^{\circ}_{1}(\delta_{p} \otimes \delta_{q})
\end{eqnarray*}
and its inverse satisfies
$$
\ttt^{-1}_{1}(a \gamma(\delta_{p}) \otimes b \gamma
(\delta_{q}))=\ttt^{-1}_{1}(a \otimes b)( \gamma \otimes \gamma
)(T^{\circ}_{1})^{-1}(\delta_{p} \otimes \delta_{q}).
$$
Similarly, the map $\ttt_2$ satisfies
$$
\ttt_{2}(a \gamma(\delta_{p}) \otimes b \gamma
(\delta_{q}))=\ttt_{2}(a \otimes b)( \gamma \otimes \gamma
)T^{\circ}_{2}(\delta_{p} \otimes \delta_{q})   \\
$$
and its inverse satisfies
$$
\ttt^{-1}_{2}(a \gamma(\delta_{p}) \otimes b \gamma
(\delta_{q}))=\ttt^{-1}_{2}(a \otimes b)( \gamma \otimes \gamma
)(T^{\circ}_{2})^{-1}(\delta_{p} \otimes \delta_{q}).
$$

For $a,b \in \at$ and $\delta_{p}, \delta_{q} \in K(G)$  the counit satisfies
\begin{eqnarray*}
\et(a \gamma(\delta_{p})b \gamma (\delta_{q}))&=& m
\ttt^{-1}_{1}(a \gamma(\delta_{p}) \otimes b \gamma (\delta_{q}))\\
&=& m \ttt^{-1}_{1}(a \otimes b)( \gamma \otimes
\gamma)((T^{\circ}_{1})^{-1}(\delta_{p} \otimes \delta_{q}))\\
&=& \et(a)b
\gamma(\varepsilon_{\circ}(\delta_{p})\delta_{q})\\
&=&\et(a)\varepsilon_{\circ}(\delta_{p})b \gamma
(\delta_{q}).
\end{eqnarray*}
Since $\varepsilon_{\circ}(f)=f(e)\;\; \forall f \in K(G) $ then
$\et( a\gamma( \delta_{p}))= \et (a)$ when $e=p$ and zero otherwise\\
Also, the counit satisfies
$$
(I \otimes \et ) \dt(a \gamma (\delta_{p}))(1 \otimes b
\gamma (\delta_{q}))= a \gamma( \delta_{p})b \gamma( \delta_{q})=ab
\gamma( \delta_{p}\delta_{q}),
$$
$$
(\et \otimes I) (a \gamma \delta_{p} \otimes 1)\dt(b
\gamma \delta_{q})=ab\gamma( \delta_{p} )\gamma( \delta_{q})=ab
\gamma( \delta_{p}\delta_{q}).
$$
 For the antipode we have
 \begin{eqnarray*}
\st(a \gamma( \delta_{p}))b \gamma( \delta_{p})&=&(\et
\otimes I)\ttt^{-1}_{1}(a \gamma(\delta_{p}) \otimes b \gamma (\delta_{q}))\\
&=&(\et \otimes I)\ttt^{-1}_{1}(a \otimes b)( \gamma \otimes
\gamma)((T^{\circ}_{1})^{-1}(\delta_{p} \otimes \delta_{q}))\\
&=&\st(a)b \gamma((\varepsilon_{\circ} \otimes
I)(T^{\circ}_{1})^{-1}(\delta_{p} \otimes \delta_{q}))\\
&=&\st(a)b \gamma(S_{\circ}(\delta_{p})\gamma( \delta_{q})\\
&=&\st(a) \gamma(S_{\circ}(\delta_{p})) b \gamma( \delta_{q}).
\end{eqnarray*}
Since $S_{\circ}(\delta_{p})=\delta_{p^{-1}}$ then  $\st(a \gamma( \delta_{p}))=\st(a)\gamma( \delta_{p^{-1}})$.\\
Also, the antipode satisfies
\begin{eqnarray*}
m(S \otimes I)\Delta(a \gamma(\delta_{p}))(1 \otimes b
\gamma(\delta_{q}))&=& \varepsilon(a \gamma(\delta_{p}))b
\gamma(\delta_{q})= \varepsilon(a) \varepsilon_{\circ}(\delta_{p})b
\gamma(\delta_{q}),\\
m(I \otimes S)(b \gamma(\delta_{q})\otimes 1) \Delta(a \gamma(\delta_{p})&=& b
\gamma(\delta_{q}) \varepsilon(a) \varepsilon_{\circ}(\delta_{p}).
\end{eqnarray*}
%===========================================================================================================
\begin{prop}
Let $\at$ be a multiplier Hopf algebra and $\gamma : K(G)
\longrightarrow M(\at)$ be a nondegenerate homomorphism from $K(G)$
into the center of $M(\at)$ satisfies
$$
\dt(\gamma(f))= (\gamma \otimes \gamma)\Delta_{\circ}(f)
$$
then this comes from a multiplier Hopf $G-$coalgebra.
\end{prop}
\begin{proof}
Define a family of algebras $A=\{A_P\}_{p\in G}$ by $A_{p}=\at \gamma(\delta_{p})$. Then $A_{p}$ is an algebra
with nondegenerate product. Define a family of homomorphism $\de=\{\Delta_{p,q}: A_{pq} \longrightarrow M(A_{p} \otimes
A_{q})\}_{p,q\in G}$ by
$$
\Delta_{p,q}(a \gamma(\delta_{pq}))=\dt(a)(\gamma(\delta_{p})
\otimes \gamma(\delta_{q})).
$$
For $a\ga(\delta_{pq})\in A_{pq}, b\ga(\delta_q)\in A_q$ we have
\begin{eqnarray*}
\Delta_{p,q}(a \gamma(\delta_{pq}))(1\otimes b\ga(\delta_q))&=& \dt(a)(1 \otimes
b)(\gamma(\delta_{p}) \otimes \gamma(\delta_{q}))\\
&=&\ttt_{1}(a \otimes b)(\gamma(\delta_{p}) \otimes
\gamma(\delta_{q})) \in A_{p} \otimes A_{q}.
\end{eqnarray*}
For $b\ga(\delta_{pq})\in A_{pq}$ and $a\ga(\delta_p)\in A_p$ we have
\begin{eqnarray*}
(a \gamma(\delta_{p}) \otimes 1)\Delta_{pq}(b \gamma(\delta_{pq}))&=&
(a \otimes 1)\dt(b)(\gamma(\delta_{p}) \otimes \gamma(\delta_{q}))\\
&=&\ttt_{2}(a \otimes b)(\gamma(\delta_{p}) \otimes
\gamma(\delta_{q})) \in A_{p} \otimes A_{q}.
\end{eqnarray*}
The coassociativity of $\Delta=\{ \Delta_{p,q}:A_{pq}
\longrightarrow M(A_{p} \otimes A_{q}) \}_{p,q\in G}$ is follows directly
because of the coassociativity of $\dt$ and that $\gamma(K(G))$
is in the center of $M(A)$. So $\de$ is a comultiplication \\
 The maps $T^{1}_{p,q}: A_{pq} \otimes A_{q}
\longrightarrow A_{p} \otimes A_{q}$ given by
\begin{eqnarray*}
T^{1}_{p,q}(a \gamma(\delta_{p}) \otimes b \gamma(\delta_{q}))&=&
\Delta_{p,q}(a \gamma(\delta_{pq})(1 \otimes
b\gamma(\delta_{q})))\\
&=& \ttt_{1}(a \otimes b)(\gamma(\delta_{p}) \otimes
\gamma(\delta_{q}))\\
&=&\ttt_{1}(a \otimes b)(\gamma \otimes
\gamma)(T^{1}_{\circ}(\delta_{pq} \otimes
\delta_{q}))
\end{eqnarray*}
are bijective with inverses
\begin{eqnarray*}
(T^{1}_{p,q})^{-1}(a \gamma(\delta_{p}) \otimes b
\gamma(\delta_{q}))&=&\ttt^{-1}_{1}(a \otimes b)( \gamma(\delta_{pq})
\otimes  \gamma(\delta_{q}))\\
&=& \ttt^{-1}_{1}(a \otimes b)( \gamma \otimes \ga)T^{1}_{\circ}(\delta_{p} \otimes
\delta_{q})
\end{eqnarray*}
Also, the maps $T^{2}_{p,q}: A_{p} \otimes A_{pq} \longrightarrow
A_{p} \otimes A_{q}$ given by
\begin{eqnarray*}
T^{2}_{p,q}(a \gamma(\delta_{p}) \otimes b \gamma(\delta_{q}))&=&
(a \gamma(\delta_{p}) \otimes 1) \Delta_{p,q}(a \gamma(\delta_{pq})) \\
&=& \ttt_{2}(a \otimes b)( \gamma \otimes \ga)T^{2}_{\circ}(\delta_{p} \otimes
\delta_{pq})
\end{eqnarray*}
are bijective with inverses
\begin{eqnarray*}
(T^{2}_{p,q})^{-1}(a \gamma(\delta_{p}) \otimes b
\gamma(\delta_{q}))&=&\ttt^{-1}_{2}(a \otimes b)( \gamma(\delta_{p})
\otimes  \gamma(\delta_{pq}))\\
&=& \ttt^{-1}_{2}(a \otimes b) ( \gamma \otimes \ga)((T^{2}_{\circ})^{-1}(\delta_{p} \otimes
\delta_{pq}))
\end{eqnarray*}
Thus $(A=\{ A_{p} \}_{p \in G} , \Delta=\{ \Delta_{p,q}
\}_{p,q \in G} )$ is a multiplier Hopf $G-$coalgebra.\\

Now let $(\ag,\dg)$ be the multiplier Hopf algebra associated to
$(A,\de)$ as in lemma$(1.4)$. Define a map $\varphi : \at \longrightarrow \ag$ by
$$
\varphi(a)(p)=a \gamma(\delta_{p}) \in A_{p}.
$$
Since $a\in\oplus_p A_p$ for finite number of $p$'s then $\varphi(a)$ will has finite support
 (i.e. $\varphi$ is well
defined ($\varphi(\at) \subset \ag$)).\\
If $\varphi(a)= \varphi(b), \;\; a,b \in \at$ then
$$
a \gamma(\delta_{p}) =b \gamma(\delta_{p}) \; \forall p
\Longrightarrow (a-b)\gamma(\delta_{p})=0 \; \forall \delta_{p}
$$
By the nondegenerancy of $\gamma$ we have $a-b=0
\Longrightarrow a=b$ hence $\varphi$ is $1-1$.\\
Also,
$$
\forall \; b \in \ag \; \exists \; a= \sum_{p \in G}b(p) \in \at \;
\textrm{such that} \; \varphi(a)=b,
$$
then $\varphi$ is onto.\\
$\varphi$ is a nondegenerate homomorphism, where
$$
\varphi(ab)(p)=ab \gamma (\delta_{p})=a\gamma (\delta_{p}) b \gamma
(\delta_{p})= (\varphi(a) \varphi(b))(p) \; \forall \; p \in G.
$$
Moreover, we also have
\begin{eqnarray*}
\Delta_{G}(\varphi(a))(p,q)&=&
\Delta_{p,q}(\varphi(a)(pq))=\Delta_{p,q}(a \gamma( \delta_{pq}
))\\
&=& \dt(a)( \gamma(\delta_{p}) \otimes
\gamma(\delta_{q}))=((\varphi \otimes \varphi) \Delta(a))(p,q)
\end{eqnarray*}
Thus $\at$ is isomorphic to $\ag$.\\
There exist a map $\gamma_G: K(G) \longrightarrow \ag$ defined by
$$
\gamma_G(f)(p)=f(p)\gamma(\delta_{p}).
$$
Since  $\gamma(\delta_p)$ is the unit of $M(A_p)$. Then $\gamma_G(f)(p)=f(p)1_p$  and hence as in lemma $(1.4)$  $\gamma_G$ is a nondegenerate homomorphism compatible with the comultiplication and $\gamma_G(K(G))$ is    in the central of $M(A)$.

\end{proof}

The last two propositions show that there exist a one to one correspondence between the class of multiplier Hopf $G-$coalgebras and the class of multiplier Hopf algebras such that each \mha endowed with a nondegenerate hommomophism form the algebra of complex valued functions with finite support on the group $G$ into the center of its multiplier algebra and compatible with its comultiplication.
%==========================================================================================================
\begin{lem}
Let $(A,\Delta)$ be a multiplier Hopf $G-$coalgebra. Then there exist a homomorphism $\varepsilon:A_{e} \longrightarrow \complex$ such that
\begin{itemize}
    \item[] $ (I \otimes \varepsilon)((a\otimes 1)\Delta_{p,e}(b))=ab,$
    \item []$ (\varepsilon \otimes I)(\Delta_{e,p}(a)(1\otimes b))=ab $
    \end{itemize}
for all $ a,b \in A_{p},\; p\in G$.
\end{lem}
\begin{proof}
Consider the algebra $(\at,\dt)$ be as above. For  $a,b\in \at$ we have
$$
\varepsilon(a)b=\sum_{p}\sum_{q}
\varepsilon(a\gamma(\delta_{p}))b\gamma(\delta_{q})=
 \varepsilon (a\gamma(\delta_{e}))b.
$$
Define a map $\eps:A_e\to \complex $ by $\eps(a)=\et(a)\;\; \forall a\in A_e$. By the properties of the map $\et$ we have that $\eps$ is a nondegenerate homomorphism. Let $x,y \in A_p$ with $x=a\ga(\del_p), y=b\ga(\del_p)$ where $a,b \in \at$ then
\begin{eqnarray*}
(I\otimes\eps )((x\otimes 1)\de_{p,e}(y))&=& (I\otimes\eps )((a\ga(\del_p)\otimes 1)\de_{p,e}(b\ga(\del_p)))\\
&=& (I\otimes\et )((a\otimes 1)(\ga\otimes\ga)(\del_p\otimes\del_e )\dt(b))\\
&=& (I\otimes\et )((a\otimes 1)(\ga\otimes\ga)((\del_p\otimes 1)\de_\circ(\del_p ))\dt(b))\\
&=& (I\otimes\et )((a\ga(\del_p)\otimes 1)\dt (b\ga(\del_p)))\\
&=& a\ga(\del_p) b\ga(\del_p) = xy.
\end{eqnarray*}
Similarly, we can prove that $ (\varepsilon \otimes I)(\Delta_{e,p}(x)(1\otimes y))=xy$.\\
Also, for the multiplier Hopf algebra $(\ag,\dg)$ we get that $\eg(a)=\eps(a(e))\;\; \forall a \in A_G$.
\end{proof}
%=========================================================================================================
\begin{lem}   \label{ae}
If $(A,\Delta)$ is a multiplier Hopf $G-$coalgebra then there exist a  family of  antihomomorphisms $S=\{S_{p}:A_{p}\longrightarrow M(A_{p^{-1}})\}_{p \in G}$ such that for all $a\in A_{e},\;\;\;b\in A_{p}$ and $\p\in\,G$ we have
\begin{itemize}
    \item[] $ m_{p}((I\otimes S_{p^{-1}})((b \otimes 1)\Delta_{p,p^{-1}}(a)))=\varepsilon(b)a$,
     \item[] $ m_{p }((S_{p ^{-1}}\otimes I)\Delta _{p^{-1},p }(a)(1\otimes b))=\varepsilon (a)b$.
\end{itemize}
\end{lem}
\begin{proof}
In the algebra $(\at,\dt)$
$$
\st(a\ga(\del_p))=\st(a) \ga(S_{\circ}(\del_{p})).
$$
For the multiplier Hopf $G-$coalgebra $(A,\de)$ we define
$$
S_{p}: A_{p} \longrightarrow M(A_{p^{-1}})
$$
by
$$
S(a\ga(\del_p))=S(a)\ga(S_{\circ}(\del_{p})).
$$
Let $a \in A_{e}$ $b \in A_{p}$. Then
\begin{eqnarray*}
m_{p}(I \otimes S_{p^{-1}})[(b \ga(\del_{p}) \otimes 1)\de_{p,p^{-1}}(a\ga(\del_{e}))]&=&
m(I \otimes S)[(b  \otimes 1)\dt(a)(\ga \otimes \ga)(\del_{p} \otimes 1) \de_{\circ}(\del_e)] \\
&=&m(I \otimes S)[(b \ga(\del_{p}) \otimes 1)\dt(a\ga(\del_{e}))] \\
&=& \eps(a\ga(\del_{e}))b \ga(\del_{p}).
\end{eqnarray*}
Also, we can prove that
$$
m_{p}(S_{p^{-1}} \otimes I)(\de_{p^{-1},p}(a)(1 \otimes b))=\eps(a)b
.$$
For the algebra $(A_{G},\de_{G})$, we have
$$
S(a)(p) = S_{p^{-1}}(a(p^{-1})).
$$
\end{proof}
%=====================================================================================================
\begin{thm}
 $A$ is a Hopf $G-$coalgebra if and only if $A$ is multiplier Hopf $G-$coalgebra such that $A_p$ is unital $\forall p\in G$.
\end{thm}

\begin{proof}
 If  $(A,\Delta)$ is a Hopf $G-$coalgebra then  $A_{p}$ has an identity $\forall p \in G $ which means that
$ M(A_{p})= A_{p} \;\;\; \forall p \in G $ and according to remark (1) we have $\Delta$ is a comultiplication.  \\
Now we will prove that the linear maps
\begin{eqnarray*}
&T_{p,q}^{1}:A_{pq}\otimes A_{q}\longrightarrow A_{p}\otimes A_{q}&,  \\
&T_{p,q}^{2}:A_{p}\otimes A_{pq}\longrightarrow A_{p}\otimes A_{q}&
\end{eqnarray*}
defined by
\begin{eqnarray*}
&T_{p,q}^{1}(a \otimes b)=\Delta_{p,q}(a)(1\otimes b)&,  \\
&T_{p,q}^{2}(a \otimes b)=(a\otimes 1)\Delta_{p,q}(b)&
\end{eqnarray*}
are bijective for all $p,q \in G\;.$\\
Consider the following linear maps
\begin{eqnarray*}
&R_{p,q}^{1}:A_{p}\otimes A_{q}\longrightarrow A_{pq}\otimes A_{q}&,  \\
&R_{p,q}^{2}:A_{p}\otimes A_{q}\longrightarrow A_{p}\otimes A_{pq}&
\end{eqnarray*}
which defined by
\begin{eqnarray*}
&R_{p,q}^{1}(a\otimes b)=((I\otimes S_{q^{-1}})\Delta_{p q,q^{-1}}(a))(1\otimes b)&,   \\
&R_{p,q}^{2}(a\otimes b)=(a\otimes 1)((S_{p^{-1}}\otimes I)\Delta_{p^{-1},p q}(b))&.
\end{eqnarray*}
By using the properties of $S$ and $\Delta$ we can prove that $R^{1}_{p,q}$ and  $R^{2}_{\p,\q}$ are the inverses of $T^{1}_{\p,\q}$ and $T^{2}_{\p,\q}$, respectively. By using  Sweedler's notation, we get
\begin{eqnarray*}
T^{1}_{\p,\q}R^{1}_{\p,\q}(a\otimes b)&=& \sum_{(a)}T^{1}_{\p,\q}(a_{(1,\p\q)} \otimes S_{\q^{-1}}(a_{(2,\q^{-1})})\,b)   \\
&=& \sum _{(a)}a_{(1,\p)}\otimes a_{(2,\q)} S_{\q^{-1}}(a_{(3,\q^{-1})})\,b  \\
&=& \sum _{(a)}a_{(1,\p)}\otimes m_{\p}(I \otimes S_{\q^{-1}})\Delta_{\q,\q^{-1}} (a_{(2,1)})(1\otimes b)  \\
&=& \sum _{(a)}(a_{(1,\p)}\otimes \varepsilon(a_{(2,1)}) 1_{\q})(1 \otimes b)=\,a\otimes b.
\end{eqnarray*}
Similarly, for $R_{\p,\q}^{1}T_{\p,\q}^{1},R_{\p,\q}^{2}T_{\p,\q}^{2}$ and $T_{\p,\q}^{2}R_{\p,\q}^{2}$. Therefore   $A$ is a multiplier Hopf $\pi-$coalgebra.\\

If $A$ is a  multiplier Hopf $G-$coalgebra such that $A_p$ is unital $\forall p\in G$ then $ M(A_{\p}\otimes A_{\q})=A_{\p}\otimes A_{\q} \; \forall \:\p,\q \in G $ and by using remark (1) we have
$$
(\Delta _{\p ,\q }\otimes I)\Delta _{\p \q ,r }=(I\otimes \Delta _{\q ,r })\Delta _{\p ,\q r }
$$
as a maps from $A_{\p\q r}$ into $A_{\p}\otimes A_{\q}\otimes A_{r}$. \\
Since $A$ is unital then $ M(A_{1}\otimes A_{\p})=A_{1}\otimes A_{\p}$ and $M(A_{\p}\otimes A_{1})=A_{\p}\otimes A_{1}$ so the counitary property just means
$$ (I\otimes \varepsilon)\Delta_{\p,1}= (\varepsilon \otimes 1)\Delta_{1,\p}=I\;\;\;\; \forall \p\in G. $$
Then we obtain  the definition of the counit introduced by Turaev.\\
If $A$ is unital then the antihomomorphism $S$ has the form
$$ S=\{S_{\p}:A_{\p}\longrightarrow A_{\p^{-1}}\}_{\p \in G}$$ such that
$$ m_{\p}(S_{\p^{-1}}\otimes I)\Delta_{\p^{-1},\p}=\varepsilon 1_{\p}=m_{\p}(I\otimes S_{\p^{-1}})\Delta_{\p,\p^{-1}}.$$
\end{proof}
%========================================================================================
\subsection{Regular multiplier Hopf group coalgebra.}
\begin{lem}
$(A=\{A_{p}\}_{p \in G},\Delta) $ is a regular \mhgc \textit{iff} $(A_{G},\de_{G})$ is a regular \mha.
\end{lem}
\begin{proof}
1) If $(A=\{A_{p}\}_{p \in G},\Delta)$ is a regular \mhgc, then  $(A_{G},\de_{G})$ is a \mha. It is clear that $\de_{G}$ is  nondegenerate  and from the coassociativity of $\de_{G}$, it follows that $\de_{G}'$ is coassociative. Let $a,b \in  A_{G}$ then
\begin{eqnarray*}
\de_{G}'(a)(1\otimes b)(q,p)&=&\sigma(\de_{G}(a)(b\otimes 1))(q,p)\\
&=& \sigma(\de_{G}(a)(b \otimes 1)(p,q))\\
&=& \sigma (\de_{p,q}(a(pq)(b(p) \otimes 1)))\\
&=& \de'_{p,q}(a(pq)(1 \otimes b(p))) \in A_{q} \otimes A_{p}.
\end{eqnarray*}
Similarly, we can prove that $(a\otimes 1)\dg'(b)  \in A_{q}\otimes A_{p}$. Hence, $\dg$ is an opposite  comultiplication on $A_{G}$.\\
The map $T'_{1}:A_{G}\otimes A_{G} \to A_{G} \otimes A_{G} $ defined by
\begin{eqnarray*}
T'_{1}(a\otimes b)(q,p)&=& \dg'(a)(1\otimes b)(q,p)\\
&=& \de'_{p,q}(a(pq)(1 \otimes b(p))) \\
&=& T^{'1}_{p,q}(a(pq)(1 \otimes b(p)))\in A_{q}\otimes A_{p}
\end{eqnarray*}
is bijective since the members of the family $T^{'1}=\{T^{'1}_{p,q}\}_{p,q \in  G}$ are bijective .\\
Similarly, we can prove that the map $T'_{2}$ is bijective. \\
Therefor, $(A_{G},\dg)$ is a regular \mha\nolinebreak.\\
2)If $(A_{G},\dg)$ is a regular \mha, then $(A,\de)$ is a \mhgc\nolinebreak.\\
a)Let $a\in A_{pq} , b\in A_{p}$. Then
\begin{eqnarray*}
\dg'(a)(1\otimes b)(q,p)&=&\sigma_{G}(\de(a)(b \otimes 1)(p,q))\\
&=&\sigma(\de(a)(b \otimes 1)(p,q))\\
&=& \sigma(\de(a(pq))(b(p) \otimes 1))\\
&=& \de'_{p,q}(a(pq))(1 \otimes b(p)).
\end{eqnarray*}
Since $\dg'(a)(1\otimes b)(q,p)\in A_{q} \otimes A_{p}$ then $\de'_{p,q}(a)(1\otimes b)\in A_{q} \otimes A_{p}$.\\
Similarly, we can prove that $(a\otimes 1)\de'_{p,q}(b)\in A_{q} \otimes A_{p}$ for all $a \in A_{q} ,b \in A_{pq}$\\
b) Since the map $T'_1$ is bijective then $T'_{1}(p,q)$ is bijective . Similarly, for the maps $T^{'2}_{p,q} $.\\
c) Let $a \in A_{p}, b \in A_{r q p}$ and $c \in A_{r}$. Then
\begin{eqnarray*}
&[(a\otimes 1\otimes 1)&(\dg'\otimes I)(\dg'(b)(1\otimes c))] (p,q,r)\\
&=&[ (I\otimes \sigma_{G})(\sigma_{G} \otimes I)(I\otimes \sigma_{G})((1\otimes 1\otimes a)(I\otimes\dg)(\dg(b)(c\otimes 1)))](p,q,r)\\
&=& (I\otimes \sigma)(\sigma \otimes I)(I\otimes \sigma)((1\otimes 1\otimes a)(I\otimes\dg)(\dg(b)(c\otimes 1))(r,q,p))\\
&=& (I\otimes \sigma)(\sigma \otimes I)(I\otimes \sigma) (1 \otimes 1 \otimes a(p))(I \otimes \de_{q,p})(\de_{r,qp}(b(rqp)(c(r) \otimes 1)))\\
&=& (a(p) \otimes 1\otimes 1)(\de'_{q,p}\otimes I)(\de'_{r, qp }(b(rqp))(1 \otimes c(r))).
\end{eqnarray*}
On the other hand
\begin{eqnarray*}
&(I\otimes \dg')&((a\otimes 1)\dg'(b))(1\otimes 1\otimes c)(p,q,r)\\
&=& (\sigma_{G}\otimes I)(I\otimes \sigma_{G})(\sigma_{G}\otimes I)( (\dg\otimes I)((1\otimes a)\dg(b))(c\otimes 1\otimes 1))(p,q,r)\\
&=&\sigma\otimes I)(I\otimes \sigma)(\sigma\otimes I)( (\dg\otimes I)((1\otimes a)\dg(b))(c\otimes 1\otimes 1)(r,q,p))\\
&=& (\sigma\otimes I)(I\otimes \sigma)(\sigma\otimes I)((\de_{r,q}\otimes I)((1\otimes a(p) )\de_{rq,p}(b(rq,p)))(c(r)\otimes 1\otimes 1))\\
&=&( I\otimes \de'_{r,q})((a(p)\otimes 1)\de'_{rq,p}(b(rq,p)))(1\otimes 1\otimes c(r)).
\end{eqnarray*}
Since $\dg'$ is coassociative, we have
$$
(a(p)\otimes 1\otimes 1)(\de'_{q,p}\otimes I)(\de_{r,qp}(b(r,qp))(1\otimes c(r)))=( I\otimes \de'_{r,q})((a(p)\otimes 1)\de_{rq,p}(b(rq,p)))(1\otimes 1\otimes c(r)).
$$
Then $(A,\de)$ is regular \mhgc.
\end{proof}
%======================================================================================================
\begin{cor}
Let $(A,\de)$ be a \mhgc. Then for any element $a\in A_p$ there exist some elements $e,f\in A_p$ such that $ea=a=af$. If the antipode of $(A,\de)$ is endomorphism of $A$ (i.e. $S_p(A_p)\subset A_{p^{-1}}$ for all $p\in G$) then such elements $e,f$ exist for any finitely many elements $a_1,a_2,...,a_n$ in $A_p$ such that $ea_i=a_i=a_if$ for any $i$.
\end{cor}
%=========================================================================================================
\begin{cor}
If $(A,\de)$ is a regular \mhgc then
\begin{itemize}
\item [1.] $(\varepsilon \otimes I)((1\otimes a)\Delta _{e,p}(b))=ab$,
\item [2.] $(I\otimes \varepsilon )(\Delta _{p,e}(a)(b\otimes 1))=ab$
\end{itemize}
 for all $p \in G $ and $a,b \in
A_{p}$.
\end{cor}
\begin{proof}
Since $(A,\de)$ is a regular \mhgc then $(A_{G},\dg)$ is a regular \mha.
For all $a,b \in A_{p}$ we have
\begin{eqnarray*}
a b(p) &=& (I\otimes \eg)(\dg(a)(b\otimes 1))(p)\\
&=& (I\otimes \eps)(\dg(a)(b\otimes 1)(p,e))\\
&=&(I\otimes \eps)(\de_{p,e}(b\otimes 1)).
\end{eqnarray*}
Similarly, we can prove part 2.
\end{proof}
%=============================================================================================
\begin{lem}
Let $(A,\de)$ be a regular \mhgc. Then
\begin{itemize}
\item There exist an antihomomorphism $S'=\{S'_{p}:A_{p}\longrightarrow M(A_{p^{-1}})\}_{p \in G}$ such that for all $a\in A_{1}\;\;\;b,c\in A_{p}$ and $p \in G$
\begin{itemize}
    \item[] $ m_{p}((I\otimes S'_{p^{-1}})((c\otimes1)\de'_{p,p^{-1}}(a))(1\otimes b))=c \eps (a)b$,
\item[] $ m_{p}((c\otimes 1)(S'_{p ^{-1}}\otimes I)\de'_{p^{-1},p }(a)(1\otimes b))=c\eps (a)b$.
\end{itemize}
\item   $S_{p}(A_{p })\subseteq A_{p
^{-1}}$, and $S'_{p}(A_{p})\subseteq A_{p ^{-1}}$,  $S_{p}$ is invertible for all $p \in G$ with inverse $S^{'}_{p^{-1}}$  for\ all$\;p \in G $.
\end{itemize}
\end{lem}
The proof of the first part is completely similar to lemma $(1.7)$ if we use the \mha $(A_{G},\dg')$ instead of $(A_{G},\dg)$. The $2^{nd}$ part is follows directly form the properties of $S_{G}$ and $S^{'}_{G}$.
%============================================================================================================
\begin{lem} \label{a}
\begin{enumerate}
\item  Let $a\in A_{pq}\;,b\in A_{q }\;,a_{i}\in A_{p
}$, and $b_{i}\in A_{q^{-1}}$. Then
$$
a\otimes S_{q }(b)=\sum \de _{pq ,q
^{-1}}(a_{i})(1\otimes b_{i})\Leftrightarrow \sum a_{i}\otimes S_{q^{-1}
}^{'}(b_{i})=(1\otimes b)\de_{p ,q}(a).
$$
\item  Let $a\in A_{pq}\;,b\in A_{p }\;,a_{i}\in A_{q },
$ and $\;b_{i}\in A_{p ^{-1}}$. Then
$$
S_{p ^{-1}}(b)\otimes a=\sum (b_{i}\otimes 1)\de_{p
^{-1},pq}(a_{i})\Leftrightarrow \de_{p,q
}(a)(b\otimes 1)=\sum S_{p ^{-1}}^{'}(b_{i})\otimes a_{i}.
$$
\item Let $a\in A_{p ^{-1}}\;,b\in A_{q
^{-1}p ^{-1}}\;,a_{i}\in A_{p ^{-1}}$ and $b_{i}\in A_{q}$. Then the following statements are equivalent :
\begin{enumerate}
\item  $\de_{q,q ^{-1}p ^{-1}}(a)(1\otimes b)=\sum \de
_{q ,q^{-1}p ^{-1}}(a_{i})(b_{i}\otimes 1)$,
\item  $a\otimes S_{pq }^{'}(b)=\sum (a_{i}\otimes 1)\de
_{p^{-1},pq }(b_{i})$,
\item $(1\otimes a)\de_{q^{-1},p^{-1}}(b)=\sum S_{q }(b_{i})\otimes a_{i}. $
\end{enumerate}
\end{enumerate}
\end{lem}
\begin{proof}
Since $(\ag,\dg)$ is a regular \mha then for all $a\in A_{pq}\;,b\in A_{q }\;,a_{i}\in A_{p}$, and $b_{i}\in A_{q ^{-1}}$, we have
\begin{eqnarray*}
(a\otimes S_G(b))(pq,q^{-1})=(\sum \dg(a_{i})(1\otimes b_{i}))(pq,q^{-1})\Leftrightarrow (\sum a_{i}\otimes S_G(b_{i}))(p.q)=((1\otimes b)\dg(a))(p,q),\\
a(pq)\otimes S_{q}(b(q))=\sum \Delta _{pq ,q
^{-1}}(a_{i}(p))(1\otimes b_{i}(q^{-1}))\Leftrightarrow \sum a_{i}(p)\otimes S_{q^{-1}}^{'}(b_{i}(q^{-1}))=(1\otimes b(q))\Delta _{p,q }(a(pq)).
\end{eqnarray*}
Similarly, we can prove  parts 2 and 3.
\end{proof}
%=====================================================================================================
\begin{lem}
Let $b\in A_{p }$ and $a\in A_{pq }$. Then
$$
(1\otimes S_{p }(b))\Delta _{q ^{-1},p ^{-1}}(S_{p q}(a))=(S_{q }\otimes S_{p })(\Delta _{p ,q }^{'
}(a)(1\otimes b)).
$$
\end{lem}
\begin{proof}
Let $b\in A_{p}$ and $a\in A_{pq }$. Then
\begin{eqnarray*}
(1\otimes S_G(b))\dg(S_G(a))(q^{-1},p^{-1})&=&(S_G\otimes S_G)(\Delta'_G(a)(1\otimes b))(q^{-1},p^{-1}),\\
(1\otimes S_{p}(b(p)))\Delta _{q ^{-1},p^{-1}}(S_{p q
}(a(pq)))&=&\sigma((S_{p }\otimes S_{q })(\Delta _{p,q }(a(pq))( b(p)\otimes 1)))\\
&=&(S_{q }\otimes S_{p })(\Delta _{p,q }^{'
}(a(pq))(1\otimes b(p))).
\end{eqnarray*}
\end{proof}
%======================================================================================================================
\begin{prop}
If $(A,\Delta)$ is a multiplier Hopf $G-$coalgebra. Define  the vector space $A^{\star}$ by
$$ A^{\star}=\oplus_{\p \in G } A^{\star}_{\p}, $$
where $A_{\p }^{*}$ is the vector space of linear functionals on $A_{\p }$
spanned by elements of the form $a\rightarrow f(bac)$ with $b,c\in
A_{\p }$ and $f\in A_{\p }^{\prime}$. $A^*$ can be made into associative algebra under  the product
$$
(fg)(a)=(f\otimes g)(\Delta _{\p ,\q }(a))
$$
where $f\in A_{\p }^{*}$ and $g\in A_{\q }^{*}.$
\end{prop}
\begin{proof}
 Let $f\in A^{*}_{\p}$ and $g\in A^{*}_{\q}$. Then there exist $b,c\in A_{\p}$, $d,e\in A_{\q}$,  $f^{'}\in A^{'}_{\p}$, and  $g^{'}\in A^{'}_{\q}$ such that  $f(a)=f^{\prime }(bac)$ for all $a\in A_{\p}$, and $g(a)=g^{'}(dae)$ for all $a\in A_{\q}$. Assume that $a\in A_{\p,\q}$ then
\begin{eqnarray*}
(fg)(a) &=&m_{c}(f\otimes g)(\Delta _{\p ,\q }(a)) \\
&=&(f^{\prime}\otimes g^{\prime})((b\otimes d)\Delta _{\p
,\q }(a)(c\otimes e)).
\end{eqnarray*}
Since $(b\otimes d)\Delta _{\p ,\q }(a)(c\otimes e)\in $ $A_{\p
}\otimes A_{\q }\;$then $(fg)(a)\in \complex\;$. It means that $fg$ is a linear
functional on $A_{\p \q }.$ By linearity it can be extended to all $
f,g\in A^{*}.$\\
\indent By the bijectivity of $T_{\p ,\q }^{1}\;$and$\;T_{\p ,\q
}^{2}\;$we can write
\begin{eqnarray*}
b\otimes d &=&\sum_{i=1}^{n}(p_{i}\otimes 1)\Delta _{\p ,\q }(q_{i}), \\
c\otimes e &=&\sum_{j=1}^{m} \Delta _{\p ,\q }(r_{j})(1\otimes k_{j}).
\end{eqnarray*}
Then
$$
(fg)(a)=\sum_{i,j}(f^{\prime}\otimes g^{\prime})((p_{i}\otimes 1)\Delta _{\p
,\q }(q_{i}ar_{j})(1\otimes k_{j})).
$$
Define
$$
h_{i,j}:A_{\p \q }\rightarrow \complex
$$
by
$$
h_{i,j}(a)=(f^{\prime}\otimes g^{\prime})((p_{i}\otimes 1)\Delta _{\p
,\q }(a)(1\otimes k_{j})).
$$
Then
$$
(fg)(a)=\sum_{i,j}h_{i,j}(qar).
$$
i.e. $fg\in A^{*}\;\;\; \forall f,g\in A^{\star}$
and thus $(A^{*},.)$ be an algebra.\\
Let $f_{1}\in A_{\p }^{*},f_{2}\in A_{\q }^{*}$ and $f_{3}\in
A_{r}^{*}\;$ defined by $f_{i}=f_{i}(b_{i}a_{i}c_{i})\; i=1,2,3,\;$where $a_{1},b_{1},c_{1}\in A_{\p },a_{2,}b_{2},c_{2}\in A_{\q }\,
$and$\,a_{3},b_{3},c_{3}\in A_{r }$. Then

\begin{eqnarray*}
((f_{1}f_{2})f_{3})(a) &=&(f_{1}f_{2}\otimes f_{3})((1\otimes b_{3})\Delta
_{\p \q ,r }(a)(1\otimes c_{3})) \\
&=&(f_{1}\otimes f_{2}\otimes f_{3})((b_{1}\otimes b_{2}\otimes 1)((1\otimes
b_{3})\Delta _{\p \q ,r }(a)(1\otimes c_{3}))(c_{1}\otimes
c_{2}\otimes 1) \\
&=&(f_{1}\otimes f_{2}\otimes f_{3})((1\otimes b_{2}\otimes
b_{3})((b_{1}\otimes 1\otimes 1)(\Delta _{\p ,\q }\otimes I)(\Delta
_{\p \q ,r }(a)(1\otimes c_{3}))) \\
&&(c_{1}\otimes c_{2}\otimes 1)) \\
&=&(f_{1}\otimes f_{2}\otimes f_{3})((1\otimes b_{2}\otimes b_{3})((I\otimes
\Delta _{\q ,\p })((b_{1}\otimes 1)(\Delta _{\p ,\q r
}(a))(1\otimes 1\otimes c_{3})) \\
&&(c_{1}\otimes c_{2}\otimes 1)) \\
&=&(f_{1}\otimes f_{2}\otimes f_{3})((1\otimes b_{2}\otimes b_{3})((I\otimes
\Delta _{\q ,\p })((b_{1}\otimes 1)(\Delta _{\p ,\q r
}(a))(c_{1}\otimes 1)) \\
&&(1\otimes c_{2}\otimes c_{3})) \\
&=&(f_{1}(f_{2}f_{3}))(a).
\end{eqnarray*}
Then $A^{*}$ becomes an associative algebra.
\end{proof}
\begin{rem}
By direct calculations we can prove that the two algebras $(A^*,*)$ and $(A{^*}_{G},m)$ are the same, where the latter is the dual algebra of $(\ag,\dg)$.
\end{rem}

\subsection{Invariant functionals on regular multiplier \\Hopf Group coalgebras}
\noindent \textbf{Notation} Let $(A,\de)$ be  a regular \mha,
$a\in \Apq$, $b\in A_q$, and $\omega\in A_p'$. Then we can define
a multiplier $n\in M(A_q)$ by
\begin{eqnarray*}
nb=(\om\otimes I)(\dpq(a)(1\otimes b))\\
bn=(\om\otimes I)((1\otimes b)\dpq(a)).
\end{eqnarray*}
We will write
$$
n=(\om\otimes I)\dpq(a).
$$
Similarly, we can define $(I\otimes \om')\dpq(a)\in M(A_p)$
,where $\om'\in A_q'$.
%=======================================================================
\begin{defn}
A family $\fa=\{\fa_p:A_p\to \complex\}_{p\in G}$ of linear
functionals on A is called left invariant \nolinebreak if
$$
(I\otimes \fa_q) \dpq(a)=\fpq(a)1_{p}\;\;\;\;\forall a\in \Apq,\;\;
p,q\in G.
$$
 A family $\fa=\{\fa_p\}_{p\in G} $ is said
to be nonzero if $\fa_p\neq 0$ for some $p \in G$.\\
 A family $\pa=\{\pa_p:A_{p}\to \complex\}_{p\in G} $ of linear
functionals on A is called right invariant  if
$$
(\pa_p\otimes 1) \dpq(a)=\ppq(a)1_{q}\;\;\;\;\forall a\in \Apq,
p,q\in G
$$
A family $\pa=\{\pa_p\}_{p \in G} $ is said to be nonzero if
$\pa_p\neq 0$ for some $p\in G$.
\end{defn}
%==============================================================================
\begin{lem}
If $\fa$ is a nonzero left invariant family of functionals then
$\fa_p\neq 0$ for all $p \in G$ such that $A_p\neq 0$.
\end{lem}
\begin{proof}
If \fah is a nonzero left invariant family of functionals then
$\fa_q\neq 0$ for some $q\in G$. Let $ p \in G$ such that
$A_p\neq 0$. Then $A_{qp^{-1}}\neq 0$ and we have for all
$b\in A_{qp^{-1}}$
$$
(I\otimes \fa_p) \de_{qp^{-1},p}(a)(b\otimes 1)=\fa_q(a)b.
$$
Hence $\fa_p\neq 0$.\\
Similarly, the right case can be proven.
\end{proof}
%=========================================================================
\begin{prop}
$A$ has  a left invariant family of functionals iff $\ag$ has a
left invariant functional.
\end{prop}
\begin{proof}
Let \fah be a left invariant family of functionals for $A$. Define a linear functional $\fag: \ag\to \complex $ by$$ \fag(a)=\sum_{p\in G} \fa_p(a(p))\hspace{.5in}\forall a\in \ag.$$
For $a,b \in \ag$, we have
\begin{eqnarray*}
[(I\otimes \fag)(\dg(a)(b\otimes 1))](p)&=&\sum_{q\in G}(I\otimes \fa_q)[(\dg(a)(b\otimes 1))(p,q)]\\
&=&\sum_{q\in G}(I\otimes \fa_q)[\dpq(a(pq))(b(p)\otimes 1)]\\&=& \sum_{q\in G}\fa_{pq}(a(pq))b(p)\\
&=& (\fag(a)b)(p).
\end{eqnarray*}
Similarly, we can prove  $(I\otimes \fag)((b\otimes
1)\dg(a))=b\fag(a)$. Thus $\fag$ is left invariant.\\

Let $\fag$ is a left invariant functional for $\ag$. Define
$\fa_p:A_p\to \complex $ by $\fa_p(a(p))=\fag(a)$ for all $a\in A_p$.
Let $ a\in \Apq$ and $ b\in A_p$. Then
\begin{eqnarray*}
(I\otimes \fa_q)\dpq(a(pq))(b(p)\otimes 1)&=& (I\otimes \fa_q)([\dg(a)(b\otimes 1)](p,q))\\
&=& \sum_{r\in G}(I\otimes \fa_r)([\dg(a)(b\otimes 1)](p,r))\\
&=& [(I\otimes \fag)(\dg(a)(b\otimes 1))](p)\\
&=&  \fag(a)b(q)=\fa_{pq}(a)b(q).
\end{eqnarray*}
Note $(I\otimes \fa_r)([\dg(a)(b\otimes 1)](p,r))= 0\;\; \forall r\neq q $. \\
Similarly, we can prove that $(I\otimes \fa_p)( (b\otimes
1)\de_{p,q}(a))= b\fa_{pq}(a)$. Thus $\fa$ is left invariant.\\
Similarly, the right case can be proven.
\end{proof}
%=====================================================================================
\begin{cor}
Let $(A,\de)$ be a regular \mhgc with left invariant family of functionals. Then for any finitely many elements $a_q,a_2,...,a_n\in A_p$ there exist an element $e\in A_p$ such that such that $ea_i=a_i=a_ie$ for any $i$
\end{cor}
%===================================================================================
\begin{cor}
If $\fa$ is left invariant then $\fa\circ S=\{\fa_{p^{-1}}\circ
S_{p}:A_P\to \complex\}_{\pg}$ is right invariant.
\end{cor}
\begin{proof}
Let \fah be a left invariant family of functionals for $A$. Then $\fag$ is a left invariant functional for $\ag$. Hence, $\fag\circ S_G$ is a right invariant functional for $A_G$. Then the family $\{\fa_{p^{-1}}\circ S_{p}:A_p\to\complex\}_{\pg}$ is right invariant.
\end{proof}
%=======================================================================================
\begin{cor}
If \fah is a nonzero left invariant family of functionals for $A$
then $\fa$ is unique (up to scaler). The members of $\fa$ are
faithful. Also, if \pah is a nonzero right invariant family of
functionals for $A$ then $\pa$ is unique (up to scaler). The
members of $\pa$ are faithful.
\end{cor}
The proof of uniqueness is an immediate consequence of the
uniqueness of  invariant functionals on $\ag$. If $a\in A_p$ such
that $\fa_p(a(p)b(p))=0 \;\;\forall b\in A_p$ then  $\fag(ax)=0
\;\;\forall\; x\in \ag$ and hence $a=0$ ($\fag$ is faithful ). Similarly,
If $b\in A_p$
such that $\fa_p(a(p)b(p))=0 \;\;\forall \;\;a\in A_p$ then $b=0$.\\
If $A_p$ is unital for all $\pg$ then $A$ is a Hopf $G-$coalgebra and the above corollary gives the uniqueness of its left integrals.\\

Let $(A,\de)$ be \mhgc with a left invariant family of functionals  \fah. Then $(\ag,\dg)$ has a left invariant functional $\fag$. Hence, there exist a multiplier $\me_G\in M(\ag)=\{a:G\to \cup M(A_p)\mid a(p)\in M(A_p)\}$ such that
$$
(\fag\otimes I)\dg(a)=\fag(a)\me_G.
$$
Let $a\in \Apq$ and $ b\in A_q$. Then
$$
[(\fag\otimes I)(\dg(a)(1\otimes b))](q)=\fag(a)\me_G(q)b(q).
$$
By the definition of $\fag$ we get
\begin{eqnarray*}
\fa_{pq}(a(pq))\me_G(q)b(q) &=&\sum_{r\in G} (\fa_{r}\otimes I)(\de_{r,q}(a(rq))(1\otimes b(q))) \\
&=&(\fa_{p}\otimes I)(\de_{p,q}(a(pq))(1\otimes b(q))).
\end{eqnarray*}
We define a family of multipliers $\me=\{\me_p\in M(A_{p})\}_{\pg}$ by $\me_p=\me_G(p)$.\\
\begin{lem}
For all $p\in G$ there exist a multiplier $\me_p\in
M(A_p)$ such that for all $p,q\in G$
$$
(\fa_p\otimes I)\de_{p,q}(a)=\fa_{pq}(a)\me_q.
$$
Also, $\me_p$ is invertible and satisfies
$$\de_{p,q}(\me_{pq})=\me_p\otimes\me_q, \hspace{.4in}\eps(\me_1)=1,\hspace{.4in}S_p(\me_p)
=\me_{p^{-1}}^{-1}.$$
\end{lem}
\begin{cor}\label{3.8}
Let \fah be a nonzero left invariant family of functionals for
$A$. Then
\begin{enumerate}
    \item For all $\pg$, $a\in A_p$
$$\fa_{p^{-1}}(S_p(a))=\fa_p(a\me_p).$$
    \item For all $a\in A_p$ there exist $b\in A_p$ such that
    $$
\fa_p(ca)=\fa_p(bc) \;\;\;\;\;\;\;\; \forall c\in A_p.
    $$
\end{enumerate}
\end{cor}
\begin{rem} Let \fah be a nonzero left invariant family of functionals for
$A$. Then
\begin{enumerate}
    \item If $\pa=\fa\circ S$, and if we replace $a$ by $ca$, $c\in A_p$ in
    part (1) of corollary (\ref{3.8}) , we have  $\pa_p(ca)=\fa_p(ca\me_p)$, with
    $b=a\me_p$ we get $\pa_p(ca)=\fa_p(cb)$.
    \item In part (2)  of corollary (\ref{3.8}) if we apply $S$ we get $\pa_p(ca)=\pa_p(bc)$ for all
    $c\in A_p$.
    \item By  (1),(2) we get that for all $a\in A_p$
    there exist $b\in A_p$ such that $ \fa_p(ac)=\fa_p(cb)$ for all
    $c\in A_p$.
    \item If we combine all of this we get the equality of the
    four following set of functionals:
    \begin{eqnarray*}
    \{\fa_p(a.)\mid a\in A_p\}, \hspace{.5in} \{\fa_p(.a)\mid a\in A_p\},\\
    \{\pa_p(a.)\mid a\in A_p\}, \hspace{.5in} \{\pa_p(.a)\mid a\in A_p\}.
    \end{eqnarray*}
\end{enumerate}
\end{rem}

Let $(A,\de)$ be \mhgc with a left invariant family of functionals \fah. Then $(\ag,\dg)$ has a left invariant functional $\fag$. Hence, there exist an automorphism $\sigma_G$  such that $\fag$ is $\sg$ invariant and $\fag(ab)=\fag(b\sg(a))$ or all $a,b\in \ag$.\\
For $a,b\in A_p$  we have
$$
\fag(ab)=\sum_{r}\fa_r(a(r)b(r))=\fa_p(a(p)b(p)).
$$
On the other hand
$$
\fag(b\sg(a))=\sum_{r}\fa_r(b(r)\sg(a)(r) )=\fa_p(b(p)\sg(a)(p)).
$$
Define $\sigma_p:A_p\to A_p$ by $\sigma_p(a(p))=\sg(a)(p)\;\; \forall\;a\in A_p$\\
Note : By the faithfulness of $\varphi_p$ we have that  $ \sigma_p(a(p))\neq 0$ for all $0\neq a\in  A_p$.
\begin{lem}
There exist a family of automorphisms $\sigma=\{\sigma_p\}_{p\in G}$ of $A$ such that $\fa_p(ab)=\fa_p(b\sigma_p(a))$ for all $a,b \in A_p$. We have also that $\fa_p$ is $\sigma_p$ invariant.
\end{lem}
\begin{cor}
Let $\sigma, \sigma'$ denoted the two families  of automorphisms associated with $\fa$ and $\pa$ respectively. Then for all $ a\in A_{pq}$ the automorphisms $\sigma, \sigma'$ satisfying :
\begin{eqnarray*}
S_p\circ \sigma'_p&=&\sigma_{p^{-1}}^{-1},\\
\dpq(\sigma_{pq}(a))&=&(S_{p^{-1}}\circ S_p\otimes \sigma_{pq})\dpq(a),\\
\dpq(\sigma'_{pq}(a))&=&(\sigma'_{pq}\otimes S_{p}^{-1}\circ S_{p^{-1}}^{-1} )\dpq(a).
\end{eqnarray*}
\end{cor}
\begin{cor}
For all $a,b\in A_p$ we have
$$
\fa_p((S_{p^{-1}}\circ S_p)(a)b)=\fa_p(b\sigma_p((S_{p^{-1}}\circ S_p)(a))).
$$
\end{cor}
\begin{proof}
let $a,b\in A_p$. Then
\begin{eqnarray*}
\fag(S^2(a)b)&=&\fa_p(((S\circ S)(a))(p)b(p)) \\
&=&\fa_p(((S_{p^{-1}}\circ S_p)(a(p)))b(p)).
\end{eqnarray*}
On the other hand, we have
\begin{eqnarray*}
\fag(b\sigma(S^2(a)))&=& \fa_p(b(p)\sigma(((S\circ S)(a)))(p)) \\
&=&\fa_p(b(p)\sigma_p((S_{p^{-1}}\circ S_p)(a(p)))).
\end{eqnarray*}
Since
$$
\fag(S^2(a)b)=\fag(b\sigma(S^2(a)))
$$
we get the required formula.
\end{proof}
%=================================================================================
Let $a\in A_p$. Then by corollary (\ref{3.8}) we have
\begin{eqnarray*}
\fa_p(S_{p^{-1}}\circ S_p(a))&=&\fa_{p^{-1}}(S_p(a)\me_{p^{-1}})\\
&=& \fa_{p^{-1}}(S_p(\me_p^{-1}a))\\
&=& \fa_{p^{-1}}(\me_p^{-1}a\me_p).
\end{eqnarray*}
By the uniqueness of the family $\fa$ then there exist a number $\tau_p$ such that $\fa_p\circ S_{p^{-1}}\circ S_p= \tau_p\fa_p$ for all $p\in G$. On the other hand  $\fag\circ S^2=\sum_p \fa_p\circ S_{p^{-1}}\circ S_p$ is left integral for $\ag$. Then $\fag\circ S^2=\tau_G\fa_g$ and hence $\sum_p \tau_p\fa_p=\sum_p\tau_G \fa_p$ which means $\tau_p=\tau_G$ for all $\pg$.
\begin{cor}
We have $\sigma_P(\me_P)=\sigma'_{p}(\me_p)=(1/\tau)\me_p$ and $\me_p\sigma(a)=\sigma'_p(a)$ for all $a\in A_p, p \in\nolinebreak G.$
\end{cor}
%=======================================================================================================
\begin{lem}
Let $(A,\de)$ be a regular \mhgc, with nonzero left invariant
family of functionals \fah. Then the vector space
$$
\Ah=\oplus_{\pg} \Ah_p,
$$
where $\Ah_p=\{\fa_p(.a)\mid a\in A_p\}$, can be made an
associative algebra under the product
$$
fg(a)=(f\otimes g)\de_{p,q}(a)
$$
for all $f\in\Ah_p, g\in \Ah_q$, and $a\in A_{p,q}$. Moreover, the
comultiplication $\deh$ given by
\begin{eqnarray*}
((f\otimes 1) \deh(g))(a\otimes b)&=& (f\otimes
g)(\de_{p,q}(a)(1\otimes b)),\\
(\deh(f)(1\otimes g))(a\otimes b) &=& (f\otimes g)((a\otimes
1)\de_{p,q}(b)),
\end{eqnarray*}
turn it into a regular multiplier Hopf algebra with nontrivial
left and right invariant functionals given by
\begin{eqnarray*}
\hat{\fa}(f)=\eps (a),&\textrm{when}& f=\pa_{\pi}(a.)\\
\hat{\pa}(f)=\eps (a),
&\textrm{when}& f=\fa_{\pi}(.a)\;.
\end{eqnarray*}
\end{lem}
\noindent The proof is direct because the multiplier Hopf algebras
$(\Ah,\deh) $ and $(\Ah_{G},\deh_{G})$ are the same
\begin{cor}
Let $(A,\de)$ be a \mhgc$\;$ with non zero left invariant family
of functionals. Then the  regular \mhgc$\;$ associated  with $(\hat{\hat{A}},\hat{\hat{\de}})$ is  isomorphic to $(A,\de)$
\end{cor}
\noindent The proof is direct since $(\ag,\dg)$ is isomorphic to
$(\hat{\hat{A_{G}}},\hat{\hat{\dg}})\cong
(\hat{\hat{A}},\hat{\hat{\de}})$ and combine this with  Proposition (1.5).

\subsection{Multiplier Hopf group coalgebras  modules}
\begin{defn} Let $A$ be a regular multiplier Hopf $G-$coalgebra.
A left $A$-module is a family of vector spaces $R=\{R_p\}_{p\in
G}$ endowed with a family of linear maps $\mu=\{\mu_p:A_p\otimes
R_p\to R_p\}_{p\in G}$ such that
$$
\mu_p(aa'\otimes x)=\mu_p(a\otimes \mu_p(a'\otimes x))
$$
for all $a,a'\in A_p$, $x\in R_p$. We call $R$ unital if
$A_pR_p=R_p$ for all $p\in G$.
\end{defn}
\begin{lem}
If $R$ is a unital left $A-$module then $R_G=\{r:G\to \cup_{p\in
G} R_P\mid r(p)\in R_p \;\;\; \forall p\in G\}$ is a unital left
$\ag-$module.
\end{lem}
\begin{proof}
Let $a\in \ag$, $r\in \rg$. We define the action of $\ag$ on $\rg$
by
$$
(ar)(p)=a(p)r(p) \;\;\;\;\;\;\;\;\; \forall p\in G.
$$
For $a,a'\in \ag, r\in \rg$, we have
$$
((aa')r)(p)=(a(p)a'(p))r(p)=a(p)(a'(p)r(p))=a(p)(a'r)(p)=(a(ar))(p).
$$
Since $A_pR_p=R_p$ then $\ag\rg=\rg$ and hence $\rg$ is a nuital
left $\ag-$module.
\end{proof}
Let $\rt$ be a unital left $\at-$module.
Then we can extend the action of $\at$ to $M(\at)$ and hence we
can define an action of $K(G)$ on $\rt$ by
$$
\del_p r= \ga(\del_p)r.
$$
Let $a \in \at, r \in \rt $ and $\del_p , \del_q \in K(G)$ then
$$
(\del_p \del_q)(a.r)= \ga(\del_p) \ga(\del_q)(a.r)=
\ga(\del_p\del_p)(a.r)=\{^{0 \;\;\;\;\;\;\;\;\;\;\;\;\;\;\;\;\; p
\neq q}_{( \ga(\del_p)a).r);\;\; p=q}\;\;.
$$
On the other hand
$$
\del_p (\del_q(a.r))= \del_p( \ga(\del_q)(a.r))=
\ga(\del_p)\ga(\del_q)(a.r)=\{^{0
\;\;\;\;\;\;\;\;\;\;\;\;\;\;\;\;\; p \neq q}_{
(\ga(\del_p)a).r\;\;\; p=q}
$$
then $(\del_p \del_q).r=\del_p (\del_q.r)$.\\
Moreover, since $\rt$ is unital $\at-$module then $\at \rt = \rt$
and hence
$$K(G). \rt =\ga(K(G)). \at \rt = \at \rt =\rt.$$
Thus $\rt$ is a unital $K(G)-$module and then $R$ can be
identified with the direct of the subspaces $\{ R_p=\del_p.R \}_{p
\in G}$. Hence, we get $a \ga(\del_p).r \in \rt_p$ $\forall r \in
\rt$
\begin{lem}
Let $\rt$ be a unital left $\at-$module. Then $\rt$ comes from a unital left $A-$module.
\end{lem}
\begin{proof}
Let $\rt$ be a unital left $\at-$module. Define  $R= \{ R_p \}_{p
\in G}$, where $R_p=\ga(\del_p)R$ and
$$
\mu_p :A_p \otimes A_p \longrightarrow R_p
$$
by
$$
\mu_p(a \ga(\del_p)\ga(\del_p)r)=\tilde{\mu}(a
\ga(\del_p)\ga(\del_p)r)
$$
For $a \ga(\del_p), b \ga(\del_p) \in A_p$, $ \ga(\del_p).r \in
R_p$ we have
\begin{eqnarray*}
\mu_p(a\ga(\del_p) b \ga(\del_p)\otimes  \ga(\del_p).r)
&=&\mt(a\ga(\del_p) b \ga(\del_p)\otimes  \ga(\del_p).r)\\
&=& \mt(a\ga(\del_p)\otimes \mt( b \ga(\del_p)\otimes  \ga(\del_p).r))\\
&=& \mt(a\ga(\del_p)\otimes \mu_p( b \ga(\del_p)\otimes
\ga(\del_p).r))\\
&=& \mu_p(a\ga(\del_p)\otimes \mu_p( b \ga(\del_p)\otimes
\ga(\del_p).r))
\end{eqnarray*}
Also, we have
$A_pR_p=\at\gamma(\del_p).\gamma(\del_p)\rt=\gamma(\del_p)\at\rt=\gamma(\del_p)\rt=R_p$.
Hence, $R$ is a unital left $A-$module\\
Define $f:\rt\to\rg$ by $f(r)(p)=\gamma(\del_p).r$ The map $f$ is
bijective moreover, it is satisfies
$$
(\varphi(a)f(r))(p)=a \ga (\del_p).r=(a.r)(p)  \;\; \forall a \in \at
$$
\end{proof}

\begin{lem}
Let $R$ be a unital left $A-$module. Then there is a unique extension
to a left $M(A)-$module. We have $1_p x=x \;\; \forall x \in
R_p,\; p \in G$.
\end{lem}
Let $R,T$ are unital left $A-$modules. Then $R \otimes T$ can be
made into a unital left $A-$module with structure map
$$
a(x \otimes y)=\de_{p,q}(a)(x \otimes y), \;\;\forall a \in A_{pq}
; \;x \in A_p;\; y \in A_q
$$
since
$$
aa'(x \otimes y)=\de_{p,q}(aa')(x \otimes y)=\de(a)\de(a')(x
\otimes y)=\de(a)(\de(a')(x \otimes y))=a(a'(x \otimes y)).
$$
Because $\de_{p,q}(A_{pq})(A_p \otimes A_q)=A_p \otimes A_q$, we
get
\begin{eqnarray*}
\de_{p,q}(A_{pq})(R_p \otimes T_q)&=&\de_{p,q}(A_{pq})(A_p \otimes
A_q)(R_p \otimes T_q)\\
&=&(A_p \otimes A_q)(R_p \otimes T_q)\\
&=& (R_p \otimes T_q).
\end{eqnarray*}
Then $R \otimes T$ is unital left $A-$module\\
Let $(A,\de)$ be a regular \mhgc, $R'$ be a $G-$graded algebra over
$\mathbb{C}$ with or without identity but with a nondegenerate
product. Assume that $R=\{R_p\}_{p \in G}$, where $R'=\oplus_{\pg} R_p$ is
a left $A-$module.
\begin{defn}
We say that $R$ is a left $A-$module algebra if
$$
a(xx')= \sum(a_{(1,p)x})(a_{(2,q)x'}) \;\; \forall a \in A, \; x
\in R_p \;, x' \in R_q.
$$
\end{defn}
\begin{lem}
If $R$ is a unital left $A-$module algebra then $R_G$ is a unital  left
$A_G-$module algebra.
\end{lem}
\begin{proof}
Let $a \in A_G$ and $x,x'\in R_G$ then
\begin{eqnarray*}
[(a_1 x)(a_2 x')]&=& \sum_r (a_1 x)_r(a_2 x')_{r^{-1}p} \\
&=& \sum_r (a_1(r) x(r))(a_2(r^{-1}p) x'(r^{-1}p)) \\
&=& \sum_r (a_{(1,r)}(p) x(r))(a_{(2,r^{-1}p)}(p) x'(r^{-1}p)) \\
&=& \sum_r a(p) (x(r)x'(r^{-1}p)) \\
&=& a(p) (xx')(p)\\
&=&[a(xx')](p).
\end{eqnarray*}
\end{proof}
Let $\rt$ be a unital left $A-$module algebra. Then $\rt$ is a
unital left $K(G)-$module algebra and hence $\rt$ is a $G-$graded
algebra.
\begin{lem}
Let $\rt$ be a unital left $\at-$module algebra. Then $\rt$ comes from a unital $A-$module algebra.
\end{lem}
\begin{proof}
Let $\rt$ be a unital left $\at-$module. Then $a\ga(\del_{pq})\in A_{pq}$, $x\ga(\del_{p}) \in A_{p}$ and $x'\ga(\del_{q}) \in A_q$
\begin{eqnarray*}
a\ga(\del_{pq})(\ga(\del_{p})x\ga(\del_{q})x')&=& \sum_r
(a_1\ga(\del_{r})\ga(\del_{p})x)(a_2\ga(\del_{r^{-1}pq})\ga(\del_{q})x') \\
&=& (\ga(\del_{p})a_1x)(\ga(\del_{q})a_2x') \\
&=&m_R(\de_{p,q}(a\ga(\del_{pq}))(\ga(\del_{p})x \otimes
\ga(\del_{q})x')).
\end{eqnarray*}
Then $R= \{R_p=\ga(\del_p)\rt \}_{p \in G}$ is a unital left $A-$module algebra.\\
Moreover, the map $f: \rt \to R_G$ defined by $f(r)(p)=
\ga(\del_p)r$ is an algebra homomorphism since
$$
f(r r')(p)=(r r')(p)= \sum_s r(s)r'(s^{-1}p)= \sum
f(r)(s)f(r')(s^{-1}p)=(f(r)f(r'))(p).
$$
\end{proof}
Let $a \in A_G$ and $x,x' \in R_G$ then we have
$$
(ax)x'= \sum a_1(x(S(a_2)x')).
$$
Then for all $p \in G$, we have
\begin{eqnarray*}
((ax)x')(p)&=&( \sum a_1(x(S(a_2)x')))(p), \\
\sum_r(ax)(r)x'(r^{-1}p)&=& \sum a_1(p)[(x(S(a_2)x'))(p)], \\
\sum_r[(a(r)x(r))x'(r^{-1}p)&=& \sum a_1(p)[\sum_r x(r)[(S(a_2)x')(r^{-1}p)] \\
&=& \sum_r \sum a_1(p)[x(r)[S(a_2)(p^{-1}r)x'(r^{-1}p)]\\
&=& \sum_r \sum a_1(p)[x(r)[S_{p^{-1}r}(a_2(p^{-1}r))x'(r^{-1}p)].
\end{eqnarray*}
Since $\dg(a)(p,q)=\de_{p,q}(a(pq))$ then $\sum a_1(p) \otimes
a_2(q)= \sum a_{(1,p)}(pq)) \otimes a_{(2,q)}(pq))$. then
$$
\sum_r(a(r)x(r))x'(p^{-1}r)= \sum_r \sum
a_{(1,p)}(r))[x(r)[S_{p^{-1}r}(a_{(2,p^{-1}r)}(r))x'(r^{-1}p)
$$
and
$$
(a(r)x(r))x'(r^{-1}p)= \sum
a_{(1,p)}(r)[x(r)[S_{p^{-1}}r(a_{(2,p^{-1}r)}(r))x'(r^{-1}p)]].
$$
\begin{lem}\begin{enumerate}
    \item For $a \in A_p$, $x\in R_p$ and $ x'\in R_q$, we have
$$
(ax)x'=\sum a_{(1,pq)}(x(S_{q^{-1}}(a_{2,q^{-1}})x')).
$$
    \item For $a \in A_q$, $x\in R_p$ and $ x'\in R_q$, we have
$$
x(ax')=\sum a_{(2,pq)}((S_{q^{-1}}(a_{1,p^{-1}})x)x').
$$
\end{enumerate}
\end{lem}

\end{document}